\documentclass{article}

\usepackage{graphicx} 
\usepackage{subfigure}

\usepackage{natbib}

\usepackage{algorithm}
\usepackage{algorithmic}
\usepackage{appendix}

\usepackage{hyperref}

\newcommand{\BEAS}{\begin{eqnarray*}}
\newcommand{\EEAS}{\end{eqnarray*}}
\newcommand{\BEA}{\begin{eqnarray}}
\newcommand{\EEA}{\end{eqnarray}}
\newcommand{\BEQ}{\begin{equation}}
\newcommand{\EEQ}{\end{equation}}
\newcommand{\BIT}{\begin{itemize}}
\newcommand{\EIT}{\end{itemize}}
\newcommand{\BNUM}{\begin{enumerate}}
\newcommand{\ENUM}{\end{enumerate}}
\newcommand{\BA}{\begin{array}}
\newcommand{\EA}{\end{array}}
\newcommand{\diag}{\mathop{\rm diag}}
\newcommand{\Diag}{\mathop{\rm Diag}}

\newcommand{\idm}{I}
\newcommand{\rb}{\mathbb{R}}
\newcommand{\R}{\mathbb{R}}
\newcommand{\BlackBox}{\rule{1.5ex}{1.5ex}}  

\newenvironment{proof}{\par\noindent{\bf Proof\ }}{\hfill\BlackBox\\[2mm]}
\newtheorem{lemma}{Lemma}

\newtheorem{proposition}{Proposition}

\newcommand{\eq}[1]{Eq.~(\ref{eq:#1})}

\input{macros}

\def \E{{\mathbb E}}
\def \P{{\mathbb P}}

\def \E{{\mathbb E}}
\def \P{{\mathbb P}}
\def \F{{\mathcal F}}

\title{A Stochastic Gradient Method with an Exponential Convergence
Rate for Finite Training Sets}

\author{
Nicolas Le Roux\\
\texttt{nicolas@le-roux.name}
\and
Mark Schmidt\\
\texttt{mark.schmidt@inria.fr}
\and
Francis Bach\\
\texttt{francis.bach@ens.fr}
\and
\\
INRIA - SIERRA Project - Team\\
D\'epartement d'Informatique de l'\'Ecole Normale Sup\'erieure\\
Paris, France
}

\date{\today}

\begin{document}
\maketitle
\begin{abstract}
We propose a new stochastic gradient method for optimizing the sum of
a finite set of smooth functions, where the sum is strongly convex.
While standard stochastic gradient methods
converge at sublinear rates for this problem, the proposed method incorporates a memory of previous gradient values in order to achieve a linear convergence
rate.  In a machine learning context, numerical experiments indicate that the new algorithm can dramatically outperform standard
algorithms, both in terms of optimizing the training error and reducing the test error quickly.
\end{abstract}

\section{Introduction}
\label{intro}

A plethora of the problems arising in machine learning involve computing an approximate minimizer of the sum of a loss function over a large number of training examples, where there is a large amount of redundancy between examples.  The most wildly successful class of algorithms for taking advantage of this type of problem structure are \emph{stochastic gradient} (SG) methods~\cite{robbins1951stochastic,bottou-lecun-2004}.
Although the theory behind SG methods allows them to be applied more generally, in the context of machine learning
SG methods are typically used to solve the problem of optimizing a sample average over a finite training set, i.e.,
\begin{equation}
\label{eq:1}
\minimize{x\in\Real^p}\quad g(x) \defd \frac{1}{n}\sum_{i=1}^n f_i(x).
\end{equation}
In this work, we focus on such \emph{finite training data} problems where
each $f_i$ is \emph{smooth} and the average function $g$ is
\emph{strongly-convex}.

As an example, in the case of $\ell_2$-regularized logistic
regression we have $f_i(x) \defd \frac{\lambda}{2}\|x\|^2 + \log(1+\exp(-b_ia_i^Tx))$, 
where $a_i \in \Real^p$ and $b_i \in \{-1,1\}$ are the training examples
associated with a binary classification problem and $\lambda$ is a
regularization parameter. More generally, any $\ell_2$-regularized empirical risk minimization problem of the
form
\begin{equation}
\label{eq:L2}
\minimize{x\in\Real^p}\quad \frac{\lambda}{2}\|x\|^2 + \frac{1}{n}\sum_{i=1}^n l_i(x),
\end{equation}
falls in the framework of~\eqref{eq:1} provided that the loss functions
$l_i$ are convex and smooth.
An extensive list of convex loss functions used in machine learning is given
by~\cite{teo2007scalable}, and  we can even include non-smooth loss functions
(or regularizers) by using smooth approximations.

The standard \emph{full gradient} (FG) method, which dates back to~\cite{cauchy1847methode}, uses iterations
of the form
\begin{equation}
\label{eq:FG}
x^{k+1} = x^k - \alpha_kg'(x^k) = x^k - \frac{\alpha_k}{n}\sum_{i=1}^nf_i'(x^k).
\end{equation}
Using $x^\ast$ to denote the unique minimizer of $g$, the FG method with a constant  step size achieves a \emph{linear} convergence rate:
\[
g(x^k) - g(x^\ast) = O(\rho^k),
\]
for some $\rho < 1$ which depends on the condition number of $g$~\cite[Theorem 2.1.15]{nesterov2004introductory}.  Linear
convergence is also known as \emph{geometric} or \emph{exponential}
convergence, because the cost is cut by a fixed fraction on each iteration.
Despite the fast convergence rate of the FG method, it can be unappealing
when $n$ is large because its iteration cost scales linearly in $n$. 
%
SG methods, on the other hand, have an iteration cost which is \emph{independent} of $n$, making them suited for that setting. The basic SG method for optimizing~\eqref{eq:1} uses iterations of the form
\begin{equation}
\label{eq:SG}
x^{k+1} = x^k - \alpha_kf_{i_k}'(x^k),
\end{equation}
where $\alpha_k$ is a step-size and a training example $i_k$ is selected
uniformly among the set $\{1,\dots,n\}$.
The randomly chosen gradient $f_{i_k}'(x^k)$ yields an unbiased estimate of
the true gradient $g'(x^k)$, and one can show under standard assumptions that,
for a suitably chosen decreasing step-size sequence $\{\alpha_k\}$, the SG
iterations achieve the sublinear convergence rate
\[
\mathbb{E}[g(x^k)] - g(x^\ast) = O(1/k),
\]
where the expectation is taken with respect to the selection of the $i_k$ variables.
Under certain assumptions this convergence rate is \emph{optimal} for strongly-convex optimization in a model of computation where the algorithm only accesses the function through unbiased measurements of its objective and gradient (see~\cite{nemirovsky1983problem,nemirovski2009robust,agarwal2010information}).
Thus, we cannot hope to obtain a better convergence rate if the algorithm only relies on unbiased gradient measurements.  Nevertheless, by using the stronger assumption that the functions are sampled from a finite dataset, in this paper we show that we can achieve an exponential converengence rate while preserving the iteration cost of SG methods.

The primay contribution of this work is the analysis of a new algorithm that we call the \emph{stochastic average
gradient} (SAG) method, a randomized variant of the incremental aggregated
gradient (IAG) method~\cite{blatt2008convergent}, which combines the low
iteration cost of SG methods with a linear convergence rate as in FG
methods. The SAG method uses iterations of the form
\begin{equation}
\label{eq:SAG}
x^{k+1} = x^k - \frac{\alpha_k}{n}\sum_{i=1}^ny^k_i,
\end{equation}
where at each iteration a random training example $i_k$ is selected and we
set
\[
y^k_i = \begin{cases}
f_i'(x^k) & \textrm{if $i = i_k$,}\\
y^{k-1}_i & \textrm{otherwise.}
\end{cases}
\]
That is, like the FG method, the step incorporates a gradient with respect
to
each training example.  But, like the SG method, each iteration only computes
the gradient with respect to a single training example and the cost of the
iterations is independent of $n$.  Despite the low cost of the SAG
iterations, in this paper we show that  \emph{the SAG iterations have a linear convergence rate}, like the FG method.  That is, by having access to
$i_k$ and by keeping a \emph{memory} of the most recent gradient value
computed for each training example $i$, this iteration achieves a faster
convergence rate than is possible for standard SG methods.
Further, in terms of effective passes through the data, we also show that for certain problems the convergence rate of SAG is faster than is possible for standard FG methods.

In a machine learning context where $g(x)$ is a \emph{training cost} associated with a predictor parameterized by $x$, we 
are often ultimately interested in the \emph{testing cost}, the expected loss on unseen data points.  
Note that a linear convergence rate for the training cost does not translate into a similar rate for the testing cost, and
an appealing
propertly of SG methods is that they achieve the optimal $O(1/k)$ rate for the \emph{testing cost} as long as every datapoint
is seen \emph{only once}.
However, as is common in machine learning, we assume that we are only given a finite training
data set and thus that datapoints are revisited multiple times.  In this context, the analysis
of SG methods only applies to the training cost and, although our analysis also 
focuses on the training cost, in our experiments
the SAG method typically reached the optimal testing cost faster than both FG and SG methods.

The next section reviews closely-related algorithms from the
literature, including previous attempts to combine the appealing aspects of
FG and SG methods.  However,  {despite $60$ years of extensive research on SG methods, most of the applications focusing on finite datasets,
we are not aware of any other SG
method that achieves a linear convergence rate while preserving the
iteration cost of standard SG methods}.  Section~\ref{convergence} states the
(standard) assumptions underlying our analysis and gives the main technical
results; we first give a slow linear convergence rate that applies for any problem, and then give 
a very fast linear convergence rate that applies when $n$ is sufficiently large.
Section~\ref{implementation} discusses practical
implementation issues, including how to reduce the storage cost from
$O(np)$ to $O(n)$ when each $f_i$ only depends on a linear combination of $x$.
 Section~\ref{experiments} presents a numerical
comparison of an implementation based on SAG to SG and FG methods, indicating that the method
may be very useful for problems where we can 
 only afford to do a few passes through
a data set.

\section{Related Work}
\label{related}


There is a large variety of approaches available to accelerate the
convergence of SG methods, and a full review of this immense literature
would be outside the scope of this work.  Below, we comment on the
relationships between the new method and several of the most closely-related
ideas.

\textbf{Momentum}: SG methods that incorporate a momentum term use
iterations of the form
\[
x^{k+1} = x^k - \alpha_k f_{i_k}'(x^k) + \beta_k(x^k - x^{k-1}),
\]
see~\cite{tseng1998incremental}.
It is common to set all $\beta_k = \beta$ for some constant $\beta$, and in
this case we can rewrite the SG with momentum method as
\[
\textstyle
x^{k+1} = x^k - \sum_{j=1}^k\alpha_j\beta^{k-j}f_{i_j}'(x^j).
\]
We can re-write the SAG updates~\eqref{eq:SAG} in a similar form as
\begin{equation}
\label{eq:SAG2}
\textstyle
x^{k+1} = x^k - \sum_{j=1}^k\alpha_kS(j,i_{1:k})f_{i_j}'(x^j),
\end{equation}
where the selection function $S(j,i_{1:k})$ is equal to $1/n$ if $j$ corresponds to the
last iteration where $j=i_k$
and is set to $0$ otherwise.
Thus, momentum uses a \emph{geometric weighting} of previous gradients while
the SAG iterations \emph{select and average} the most recent evaluation of each previous
gradient.  While momentum can lead to improved practical performance, it
still requires the use of a decreasing sequence of step sizes and is not
known to lead to a faster convergence rate.

\textbf{Gradient Averaging}:
Closely related to momentum is using the sample average of all previous
gradients,
\[
\textstyle
x^{k+1} = x^k - \frac{\alpha_k}{k}\sum_{j=1}^kf_{i_j}'(x_j),
\]
which is similar to the SAG iteration in the form~\eqref{eq:SAG} but where \emph{all} previous
gradients are used.
This approach is used in the dual averaging method~\cite{nesterov2009primal}, and while this averaging procedure 
leads to convergence for a constant step size and 
can improve the constants in the convergence rate~\cite{xiao2010dual}, it does not improve on the $O(1/k)$ rate.

\textbf{Iterate Averaging}: Rather than averaging the gradients, some authors use the basic SG iteration but take an average over $x^k$ values.
With a suitable choice of step-sizes, this gives the same asymptotic efficiency as Newton-like second-order SG methods and also leads to increased robustness of the convergence rate to the exact sequence of step sizes~\cite{polyak1992acceleration}.  Baher's method~\cite[\S1.3.4]{kushner1997stochastic} combines gradient averaging with online iterate averaging, and also displays appealing asymptotic properties.  
The epoch SG method uses averaging to obtain the $O(1/k)$ rate even for non-smooth objectives~\cite{hazan2010beyond}. 
However, the convergence rates of these averaging methods remain sublinear.

\textbf{Stochastic versions of FG methods}: Various options are available to accelerate the convergence of the FG method for smooth functions, such as the accelerated full gradient (AFG) method~\cite{nesterov1983method}, as well as classical techniques based on quadratic approximations such as non-linear conjugate gradient, quasi-Newton, and Hessian-free Newton methods.
Several authors have analyzed stochastic variants of these algorithms~\cite{schraudolph1999local,sunehag2009variable,ghadimi2010optimal,martens2010deep,xiao2010dual}. Under certain conditions these variants are convergent with an $O(1/k)$ rate~\cite{sunehag2009variable}.  Alternately, if we split the convergence rate into a deterministic and stochastic part, these methods can improve the dependency on the deterministic part~\cite{ghadimi2010optimal,xiao2010dual}.  However, as with all other methods we have discussed thus far in this section, we are not aware of any existing method of this flavor that improves on the $O(1/k)$ rate.

\textbf{Constant step size}: If the SG iterations are used with a \emph{constant} step size (rather than a decreasing sequence), then the convergence rate of the method can be split into two parts~\cite[Proposition 2.4]{nedic2001convergence}, where the first part depends on $k$ and converges linearly to $0$ and the second part is independent of $k$ but does not converge to $0$.  Thus, with a constant step size the SG iterations have a linear convergence rate up to some tolerance, and in general after this point the iterations do not make further progress.  Indeed, convergence of the basic SG method with a constant step size has only been shown under extremely strong assumptions about the relationship between the functions $f_i$~\cite{solodov1998incremental}. This contrasts with the method we present in this work which converges to the optimal solution using a constant step size \emph{and} does so with a linear rate (without additional assumptions).

\textbf{Accelerated methods}: Accelerated SG methods, which despite their name are not related to the aforementioned AFG method, take advantage of the fast convergence rate of SG methods with a constant step size.  In particular, accelerated SG methods use a constant step size by default, and only decrease the step size on iterations where the inner-product between successive gradient estimates is negative~\cite{kesten1958accelerated,delyon1993accelerated}.  This leads to convergence of the method and allows it to potentially achieve periods of linear convergence where the step size stays constant.  However, the overall convergence rate of the method remains sublinear.

\textbf{Hybrid Methods}: Some authors have proposed variants of the SG
method for problems of the form~\eqref{eq:1} that seek to gradually
transform the iterates into the FG method in order to achieve a linear
convergence rate.  Bertsekas proposes to go through the data
cyclically with a specialized weighting that allows the method to achieve a
linear convergence rate for strongly-convex quadratic functions~\cite{bertsekas1997new}.  However,
the weighting is numerically unstable and the linear convergence rate
treats full passes through the data as iterations. A related strategy is to 
 group the $f_i$ functions into `batches' of increasing size and perform SG iterations on the batches~\cite{friedlander2011hybrid}.  In both
cases, the iterations that achieve the linear rate have a cost that is not
independent of $n$, as opposed to SAG.

\textbf{Incremental Aggregated Gradient}:
Finally, Blatt et al.~presents the most closely-related
algorithm, the IAG method~\cite{blatt2008convergent}.  This method is identical to the SAG
iteration~\eqref{eq:SAG}, but uses a cyclic choice of $i_k$ rather than
sampling the $i_k$ values.  This distinction has several important
consequences.  In particular, Blatt et al.~are only able to show that the
convergence rate is linear for strongly-convex quadratic functions (without
deriving an explicit rate), and their analysis treats full passes through
the data as iterations.  
Using a non-trivial extension of their analysis 
and a proof technique
involving bounding the gradients and iterates simultaneously by a Lyapunov potential function, in this work \emph{we give an explicit
linear convergence rate for general strongly-convex
functions using the SAG iterations that only examine a single training example}.  
Further, as our analysis and experiments show, when the number of training
examples is sufficiently large, the SAG iterations achieve a linear
convergence rate under a much larger set of step sizes than the IAG method. 
This leads to more robustness to the selection of the step size and also, if
suitably chosen, leads to a faster convergence rate and improved practical
performance.  We also emphasize that in our experiments IAG and the basic FG method perform similarly, while SAG performs much better, showing that the simple change (random selection vs.~cycling) can dramatically improve optimization performance.


\section{Convergence Analysis}
\label{convergence}


In our analysis we assume that each
function $f_i$ in~\eqref{eq:1} is
differentiable and that each gradient $f'_i$ is Lipschitz-continuous with
constant $L$, meaning that for all $x$ and $y$ in $\Real^p$ we have
\[
\norm{f'_i(x)-f'_i(y)}\leq L\norm{x-y}.
\]
This is a fairly weak assumption on the $f_i$ functions, and in cases where
the $f_i$ are twice-differentiable it is equivalent to saying that the
eigenvalues of the Hessians of each $f_i$ are bounded above by~$L$.  In addition, we also assume that the average function $g =
\frac{1}{n}\sum_{i=1}^n f_i$ is strongly-convex with constant $\mu > 0$, meaning that
the function $x \mapsto g(x) - \frac{\mu}{2}\norm{x}^2$ is convex. This is a stronger assumption and is not satisfied by all machine
learning models.  However, note that in machine learning we are typically free to
choose the regularizer, and we can always add an $\ell_2$-regularization term as in Eq.~\eqref{eq:L2}
to transform any convex problem into a strongly-convex problem (in
this case we have $\mu \geq \lambda$).
Note that strong-convexity implies that the problem is solvable, meaning that there exists some unique $x^\ast$ that achieves
the optimal function value.  Our convergence results assume that we initialize $y_i^0$
to a zero vector for all $i$,
and our results depend on the variance of the gradient norms at the optimum~$x^\ast$, denoted by $\sigma^2 = \frac{1}{n} \sum_i \|f_i'(x^\ast)\|^2$.
Finally, all our convergence results consider expectations with respect to the internal randomization of the algorithm, and not with respect to the data (which are assumed to be deterministic and fixed).

We first consider the convergence rate of the method when using a constant
step size of $\alpha_k = \frac{1}{2nL}$, which is similar to the step size
needed for convergence of the IAG method in practice.
\begin{proposition}
\label{prop:small_alpha}
With a constant step size of $\alpha_k = \frac{1}{2nL}$, the SAG iterations satisfy for $k \geq 1$:
\begin{align*}
\mathbb{E}\left[\|x^k - x^\ast\|^2\right] &\leqslant \Big(1 - \frac{\mu}{8Ln}\Big)^k\Big[3\|x_0
- x^\ast\|^2 + \frac{9\sigma^2}{4L^2}\Big]  \; .
\end{align*}
\end{proposition}
The proof is given in the Appendix.  
Note that the SAG iterations also trivially obtain the $O(1/k)$
rate achieved by SG methods, since
\[
\Big(1 - \frac{\mu}{8Ln}\Big)^k \leqslant \exp\Big(- \frac{k\mu}{8Ln}\Big) \leqslant \frac{ 8 Ln}{k \mu} =  O(n/k),
\]
albeit with a constant which is proportional to $n$.  Despite this constant, they are
advantageous over SG methods in later iterations because they obtain an exponential convergence rate as in FG methods. We also note that an exponential convergence rate is obtained for any constant step size smaller than $\frac{1}{2nL}$.

In terms of passes through the data, the rate in Proposition~\ref{prop:small_alpha} is similar to that achieved by the basic FG method.
However, our
next result shows that, if the number of training examples is slightly
larger than $L/\mu$ (which will often be the case, as discussed in Section~\ref{discussion}), then the SAG
iterations can use a larger step size and obtain a better convergence rate
that is independent of $\mu$ and $L$ (see proof in the Appendix).

\begin{proposition}
\label{prop:large_alpha}
If $ n \geqslant \frac{8L}{\mu}$, with a step size of $ \alpha_k
= \frac{1}{2n\mu}$ the SAG iterations satisfy for $k \geqslant n$:
{
\begin{align*}
&\mathbb{E}\left[g(x^k) - g(x^\ast)\right] \leqslant C \Big(1 - \frac{1}{8n}\Big)^k, \\
\textrm{with } C &= \bigg[\frac{16L}{3n}\| x^0 - x^\ast\|^2 + \frac{4 \sigma^2 }{3n \mu}\Big(8 \log\Big( 1 + \frac{\mu n }{4L } \Big) + 1\Big)\bigg] \; .
\end{align*}
}
\end{proposition}
We state this result for $k \geqslant n$ because we assume that the first $n$ iterations of the algorithm use an SG method and that we initialize the subsequent SAG iterations with the average of the iterates, which leads to an $O((\log n)/k)$ rate.
In contrast, using the SAG iterations
from the beginning gives the same rate but with a constant proportional to~$n$.
Note that this bound is obtained when initializing all $y_i$ to zero after the SG phase.\footnote{While it may appear suboptimal to not use the gradients computed during the $n$ iterations of stochastic gradient descent, using them only improves the bound by a constant.} 
However, in our experiments we do not use the SG initialization but rather use a minor variant of SAG (discussed in the next section), which appears more difficult to analyze but which gives better performance. 

  It is interesting to compare this convergence rate with the known convergence rates of first-order methods~\cite[see \S2]{nesterov2004introductory}.
For example, if we take $n=100000$, $L = 100$, and $\mu=0.01$ then the basic FG method has a rate of $((L-\mu)/(L+\mu))^2 = 0.9996$ and the `optimal' AFG method has a faster rate of $(1-\sqrt{\mu/L}) = 0.9900$.  In contrast, running $n$ iterations of SAG has a much faster rate of $(1-1/8n)^n = 0.8825$ using the same number of evaluations of $f_i'$. 
Further, the lower-bound for a black-box first-order method is $((\sqrt{L}-\sqrt{\mu})/(\sqrt{L}+\sqrt{\mu}))^2 = 0.9608$, indicating that SAG can be substantially faster than any FG method that does not use the structure of the problem.\footnote{Note that $L$ in the SAG rates is based on the $f_i'$ functions, while in the FG methods it is based on $g'$ which can be much smaller.} In the Appendix, we compare Propositions~\ref{prop:small_alpha} and~\ref{prop:large_alpha} to the rates of primal and dual FG and coordinate-wise methods for the special case of $\ell_2$-regularized leasts squares.

Even though $n$ appears in the convergence rate, if we perform $n$ iterations of SAG (i.e., one effective pass through the data), the error is multiplied by $(1 - 1/8n)^n \leq \exp(-1/8)$, which is independent of $n$. Thus, each pass through the data reduces the excess cost by a constant multiplicative factor that is independent of the problem, as long as $n \geqslant 8L/\mu$.  Further, while the step size in Proposition 2 depends on $\mu$ and $n$, we can obtain the same convergence rate by using a step size as large as $\alpha_k = \frac{1}{16L}$.  This is because the proposition is true for all values of $\mu$ satisfying $\frac{\mu}{L}\geqslant \frac{8}{n}$, so we can choose the smallest possible value of $\mu = \frac{8L}{n}$.
We have observed in practice that the IAG method with a step size of $\alpha_k = \frac{1}{2n\mu}$ may diverge, even
under these assumptions.
Thus, for certain problems the SAG iterations can tolerate a much larger step
size, which leads to increased robustness to the selection of the step size.
Further, as our analysis and experiments indicate, the ability to use a
large step size leads to improved performance of the SAG iterations.

While we have stated Proposition 1 in terms of the iterates and Proposition 2 in terms of the function values, the rates obtained on iterates and function values are equivalent because,
by the Lipschitz and strong-convexity assumptions, we have $\frac{\mu}{2}\norm{x^k-x^\ast}^2 \leqslant g(x^k) - g(x^\ast) \leqslant \frac{L}{2}\norm{x^k-x^\ast}^2$.



\section{Implementation Details}
\label{implementation}
\label{sec:implementation}


In this section we describe modifications that substantially reduce the SAG iteration's memory requirements, as well as modifications that lead to better practical performance.

\textbf{Structured gradients}: For many problems the storage cost of $O(np)$ for the $y_i^k$ vectors is prohibitive, but we can often use structure in the $f_i'$ to reduce this cost. For example, many loss functions $f_i$ take the form $f_i(a_i^Tx)$ for a vector $a_i$. Since $a_i$ is constant, for these problems we only need to store the scalar $f_{i_k}'(u_i^k)$ for $u_i^k = a_{i_k}^Tx^k$ rather than the full gradient $a_i^Tf_i'(u_i^k)$, reducing the storage cost to $O(n)$.  Further, because of the simple form of the SAG updates, if $a_i$ is sparse we can use `lazy updates' in order to reduce the iteration cost from $O(p)$ down to the sparsity level of $a_i$.

\textbf{Mini-batches}: To employ vectorization and parallelism, practical SG implementations often group training examples into `mini-batches' and perform SG iterations on the mini-batches.  We can also use mini-batches within the SAG iterations, 
and for problems with dense gradients this decreases the storage requirements of the algorithm since we only need a $y_i^k$ for each mini-batch.  Thus, for example, using mini-batches of size $100$ leads to a 100-fold reduction in the storage cost.

\textbf{Step-size re-weighting}: On early iterations of the SAG algorithm, when most $y_i^k$ are set to the uninformative zero vector, rather than dividing $\alpha_k$ in~\eqref{eq:SAG} by $n$ we found it was more effective to divide by $m$, the number of unique $i_k$ values that we have sampled so far (which converges to $n$).  This modification appears more difficult to analyze, but with this modification we found that the SAG algorithm outperformed the SG/SAG hybrid algorithm analyzed in Proposition~2.

\textbf{Exact regularization}: For regularized objectives like~\eqref{eq:L2} we can use the exact gradient of the regularizer, rather than approximating it.  For example, our experiments on $\ell_2$-regularized optimization problems used the recursion
\begin{equation}
\label{eq:L2recursion}
d   \leftarrow d - y_i,  \hspace*{1cm} 
y_i  \leftarrow l_i'(x^k), \hspace*{1cm} 
d  \leftarrow d + y_i, \hspace*{1cm} 
x  \leftarrow \big(1-\alpha\lambda\big)x - \displaystyle \frac{\alpha}{m}d\; .
\end{equation}
This can be implemented efficiently for sparse data sets by using the representation $x = \kappa z$, where $\kappa$ is a scalar and $z$ is a vector, since the update based on the regularizer simply updates $\kappa$.

\textbf{Large step sizes}: Proposition~1 requires $\alpha_k \leqslant 1/2Ln$ while under an additional assumption Proposition~2 allows $\alpha_k \leqslant 1/16L$.  In practice we observed better performance using step sizes of $\alpha_k = 1/L$ and $\alpha_k = 2/(L+n\mu)$. These step sizes seem to work even when the additional assumption of Proposition~2 is not satisfied, and we conjecture that the convergence rates under these step sizes are much faster than the rate obtained in Proposition~1 for the general case.

\textbf{Line search}: Since $L$ is generally not known, we experimented with a basic line-search, where we start with an initial estimate $L_0$, and we double this estimate whenever we do not satisfy the instantiated Lipschitz inequality
\[
f_{i_k}(x^k - (1/L_k)f_{i_k}'(x^k)) \leqslant f_{i_k}(x^k) - \frac{1}{2L_k}\norm{f_{i_k}'(x^k)}^2.
\]
To avoid instability caused by comparing very small numbers, we only do this test when $\norm{f_{i_k}'(x^k)}^2 > 10^{-8}$. To allow the algorithm to potentially achieve a faster rate due to a higher degree of local smoothness, we multiply $L_k$ by $2^{(-1/n)}$ after each iteration.


\section{Experimental Results}
\label{experiments}
\label{sec:experiments}


 Our experiments compared an extensive variety of competitive FG and SG methods.  Our first experiments focus on the following methods, which
 we chose because they have no dataset-dependent tuning parameters:
\begin{list}{\labelitemi}{\leftmargin=1.7em}
   \addtolength{\itemsep}{-.215\baselineskip}
\item[--]\textbf{Steepest}: The full gradient method described by
iteration~\eqref{eq:FG}, with a line-search that uses 
cubic Hermite polynomial interpolation to find a step size satisfying
the strong Wolfe conditions, and where the parameters of the line-search
were tuned for the problems at hand.
\item [--]\textbf{AFG}: Nesterov's accelerated full gradient method~\cite{nesterov1983method}, where iterations of~\eqref{eq:FG} with a fixed step size are interleaved with an extrapolation step, and we use an adaptive line-search based on~\cite{liu2009large}.
\item[--] \textbf{L-BFGS}: A publicly-available limited-memory quasi-Newton method that has been tuned for log-linear models.\footnote{\small \url{http://www.di.ens.fr/~mschmidt/Software/minFunc.html}}  This method is by far the most complicated method we considered.
\item[--] \textbf{Pegasos}: The state-of-the-art SG method described by
iteration~\eqref{eq:SG} with a step size of $\alpha_k = 1/\mu k$ and a projection step onto a norm-ball known to contain the optimal solution~\cite{shalev2007pegasos}.
\item [--]\textbf{RDA}: The regularized dual averaging method~\cite{xiao2010dual}, another recent state-of-the-art SG method.
\item [--]\textbf{ESG}: The epoch SG method~\cite{hazan2010beyond}, which runs SG with a constant step size and averaging in a series of epochs, and is optimal for non-smooth stochastic strongly-convex optimization.
\item [--]\textbf{NOSG}: The nearly-optimal SG method~\cite{ghadimi2010optimal}, which combines ideas from SG and AFG methods to obtain a nearly-optimal dependency on a variety of problem-dependent constants.
\item[--] \textbf{SAG}: The proposed stochastic average gradient method
described by iteration~\eqref{eq:SAG} using the modifications discussed in the previous section.  We used a step-size of $\alpha_k = 2/(L_k + n\lambda)$ where $L_k$ is either set constant to the global Lipschitz constant (SAG-C) or set by adaptively estimating the constant with respect to the logistic loss function using the line-search described in the previous section (SAG-LS). The SAG-LS method was initialized with $L_0=1$ .
\end{list}

\begin{figure}
\centering

 \vspace*{-.2cm}

\hspace*{-.5cm}
\mbox{
\includegraphics[width=4.9cm]{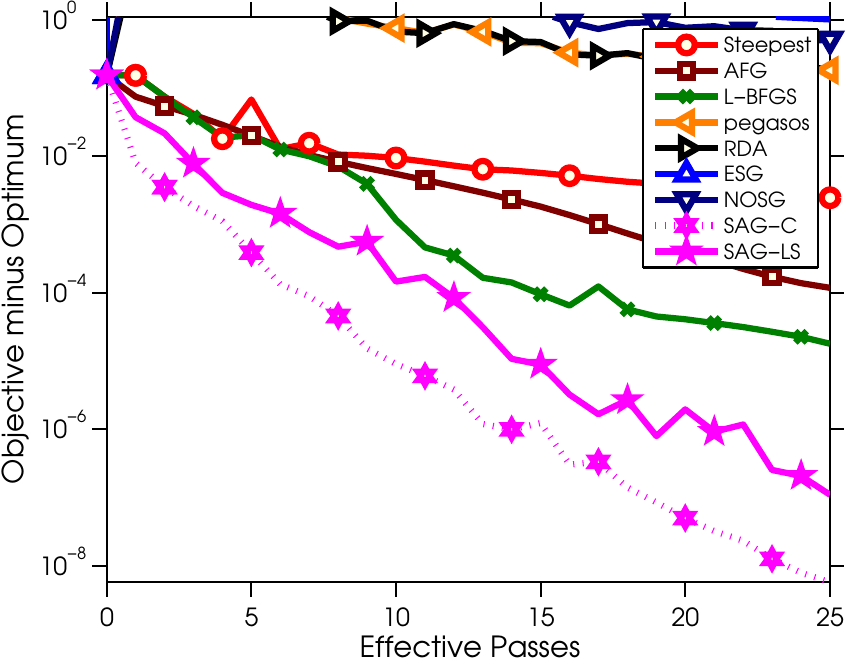} \hspace*{-.1cm}
\includegraphics[width=4.9cm]{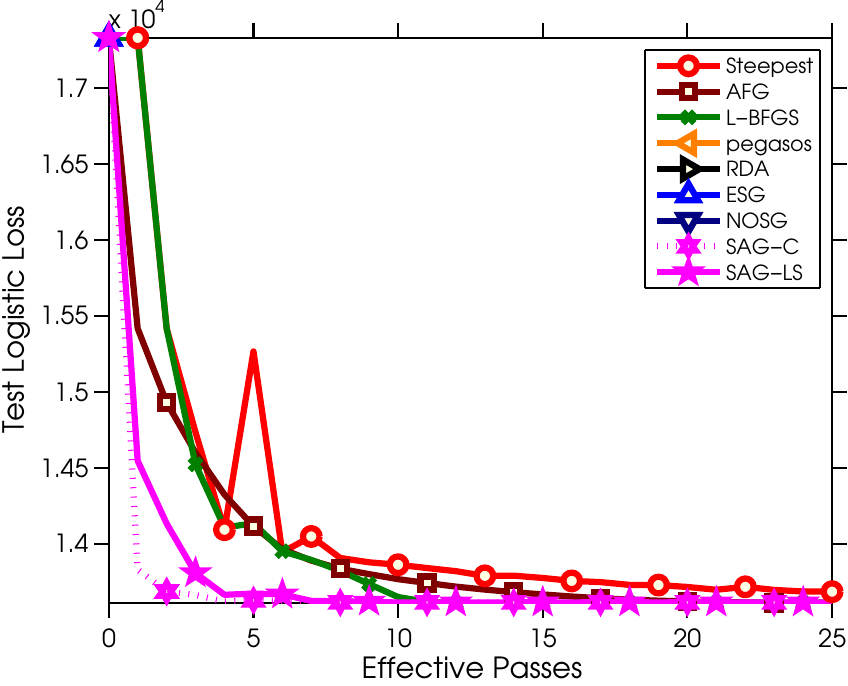} \hspace*{-.1cm}
\includegraphics[width=4.9cm]{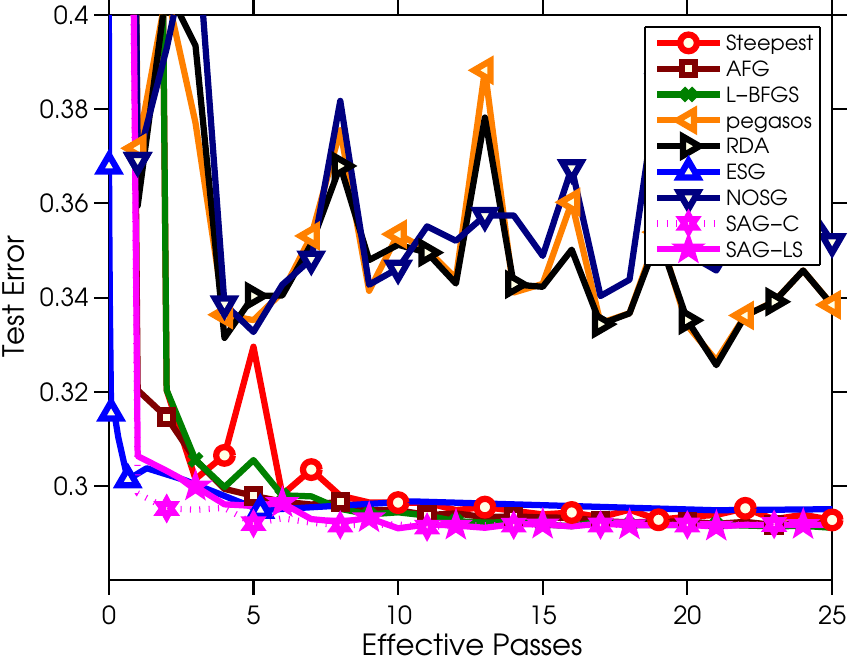}}

\hspace*{-.5cm}
\mbox{
\includegraphics[width=4.9cm]{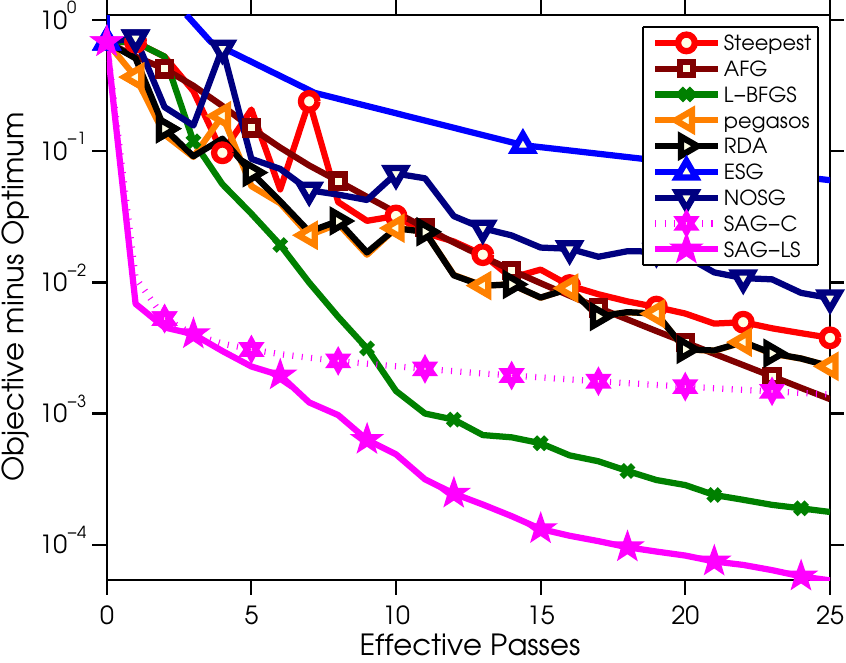} \hspace*{-.1cm}
\includegraphics[width=4.9cm]{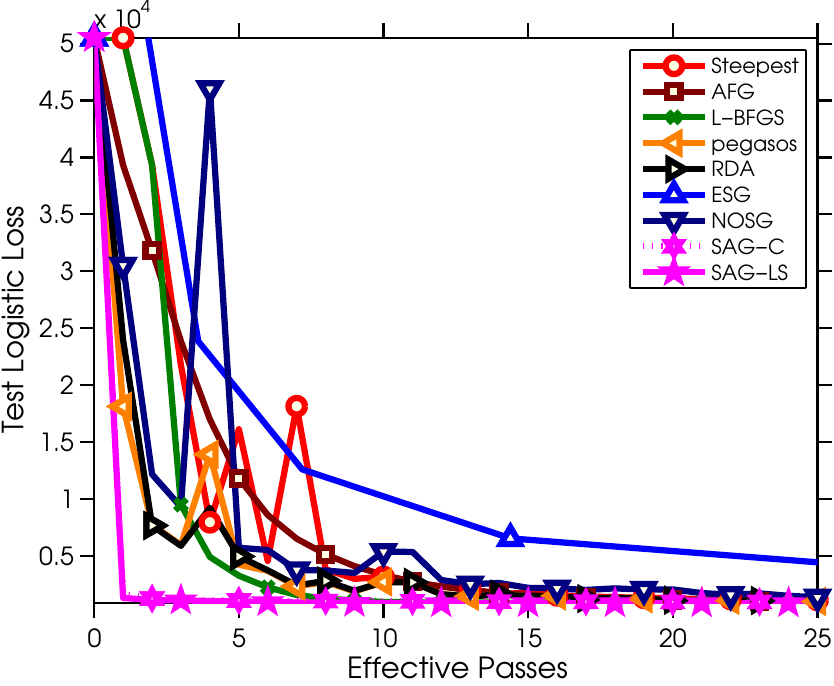} \hspace*{-.1cm}
\includegraphics[width=4.9cm]{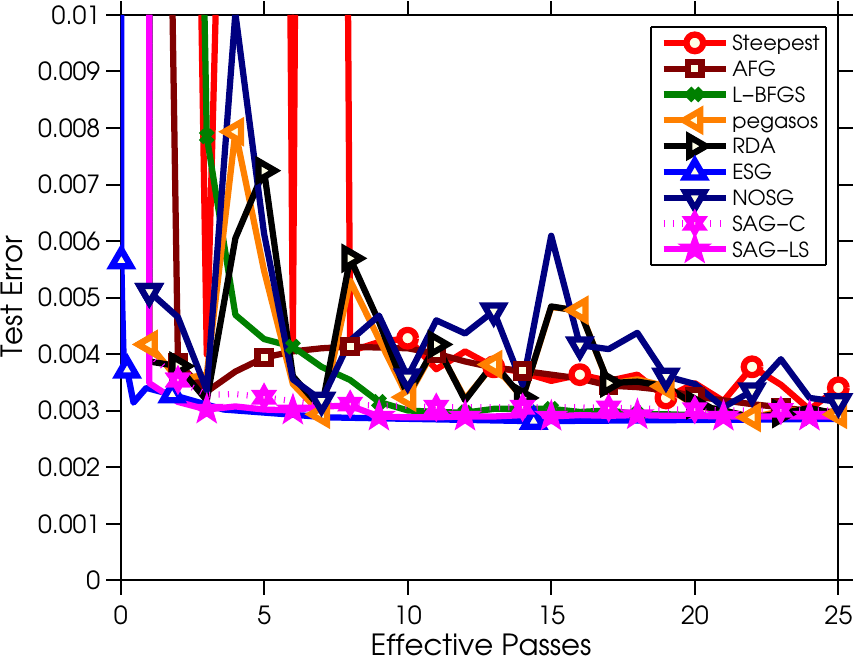}}

\hspace*{-.5cm}
\mbox{
\includegraphics[width=4.9cm]{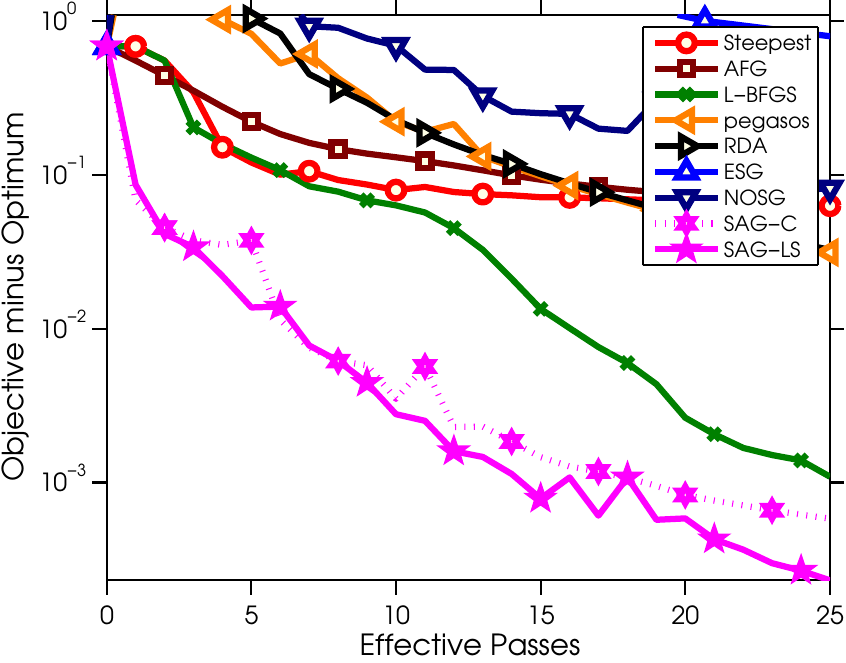} \hspace*{-.1cm}
\includegraphics[width=4.9cm]{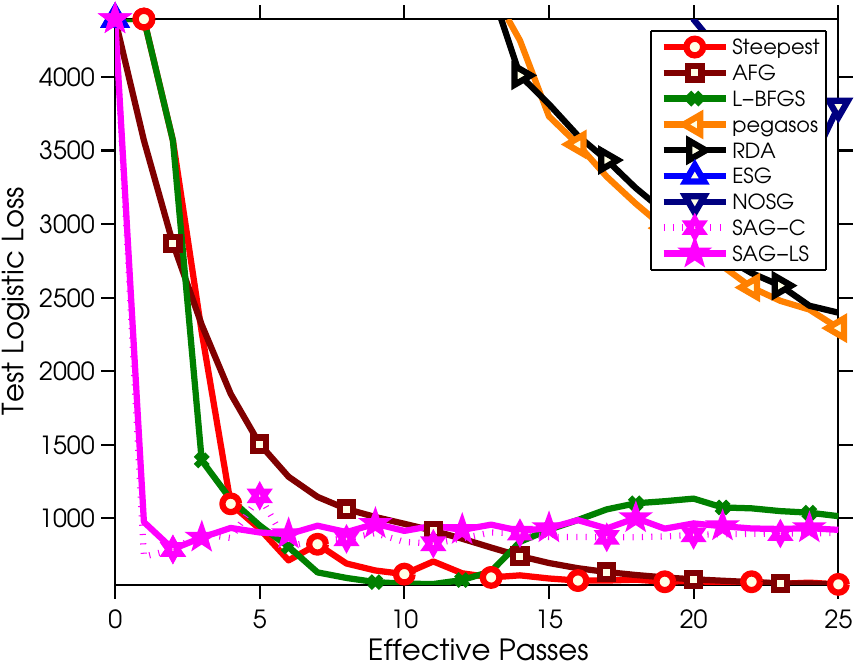} \hspace*{-.1cm}
\includegraphics[width=4.9cm]{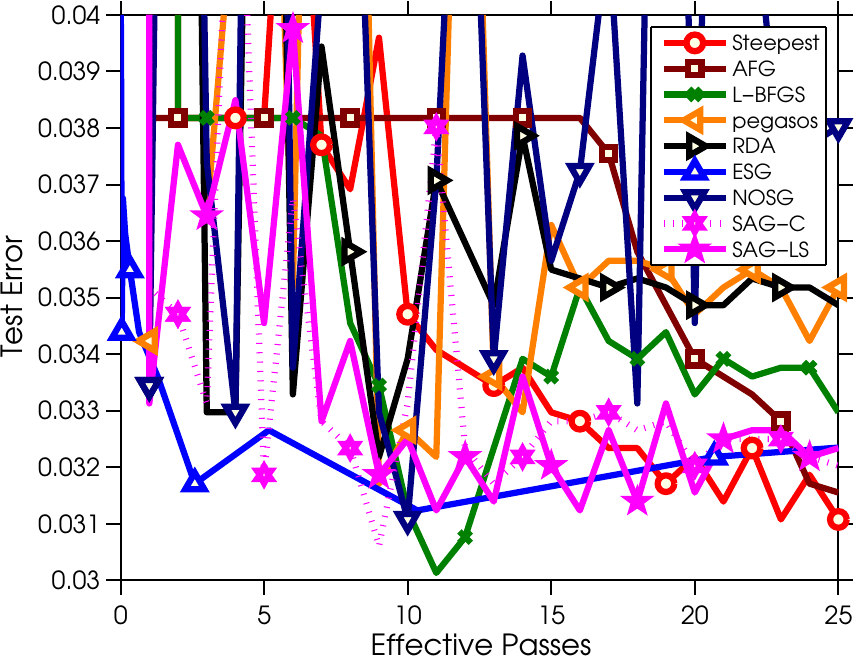}}

\hspace*{-.5cm}
\mbox{
\includegraphics[width=4.9cm]{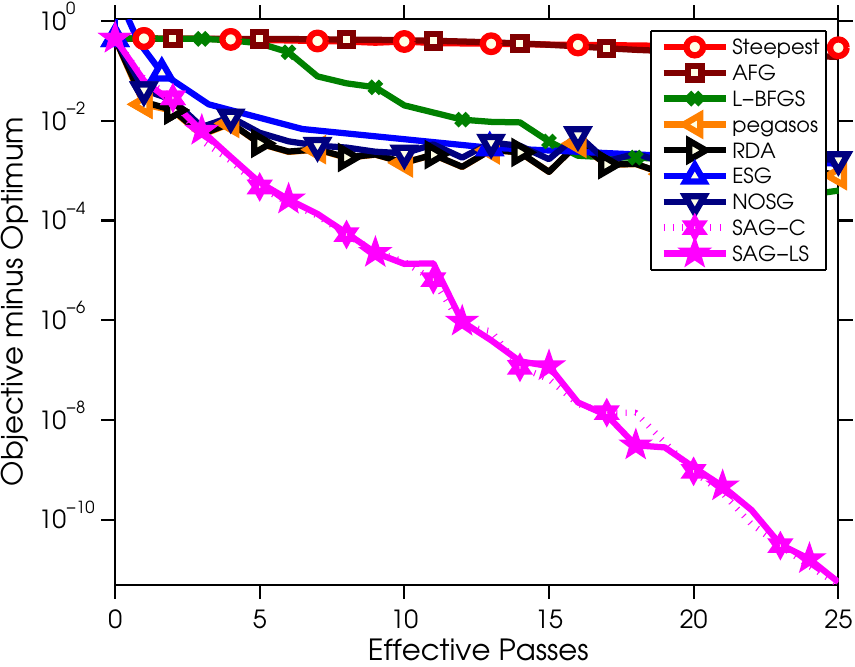} \hspace*{-.1cm}
\includegraphics[width=4.9cm]{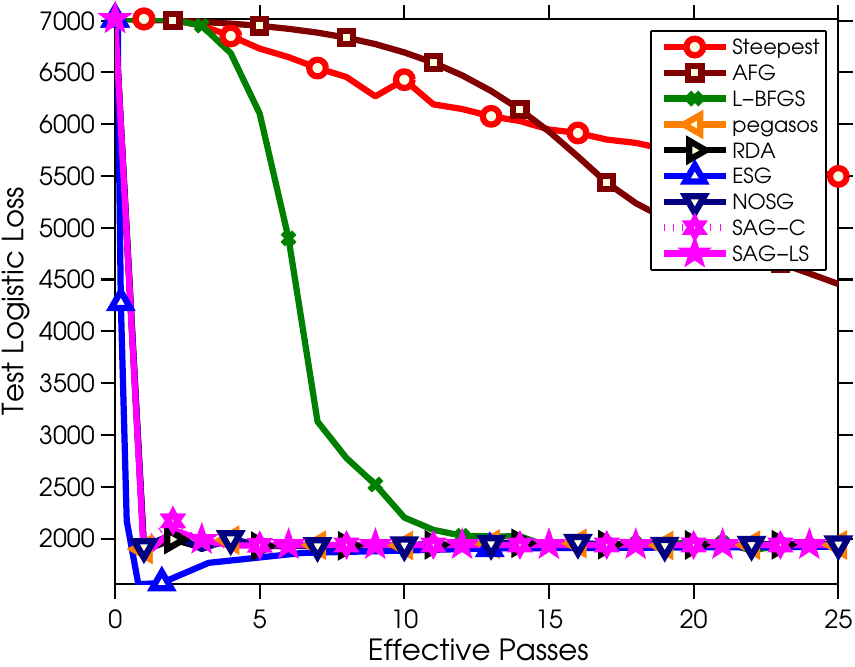} \hspace*{-.1cm}
\includegraphics[width=4.9cm]{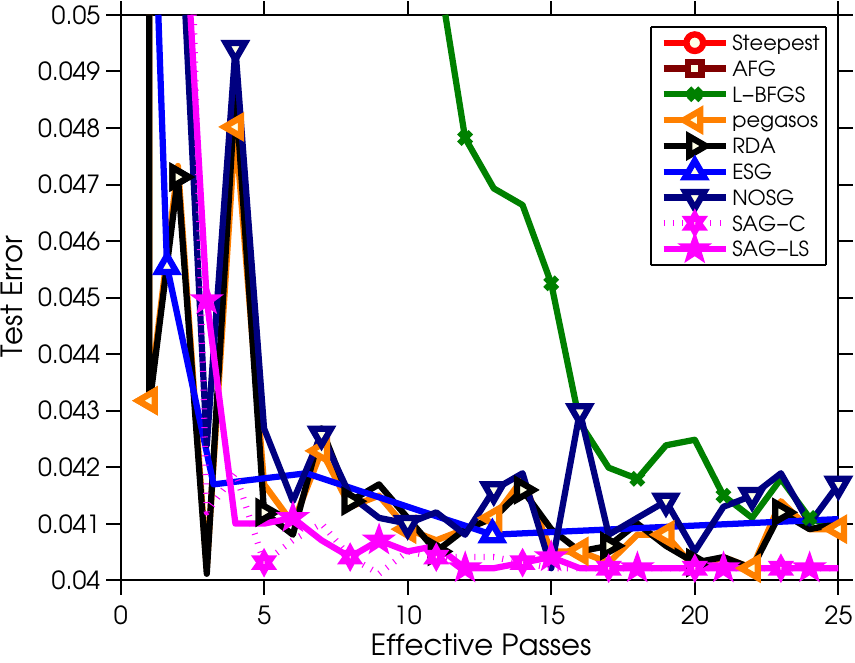}}

\hspace*{-.5cm}
\mbox{
\includegraphics[width=4.9cm]{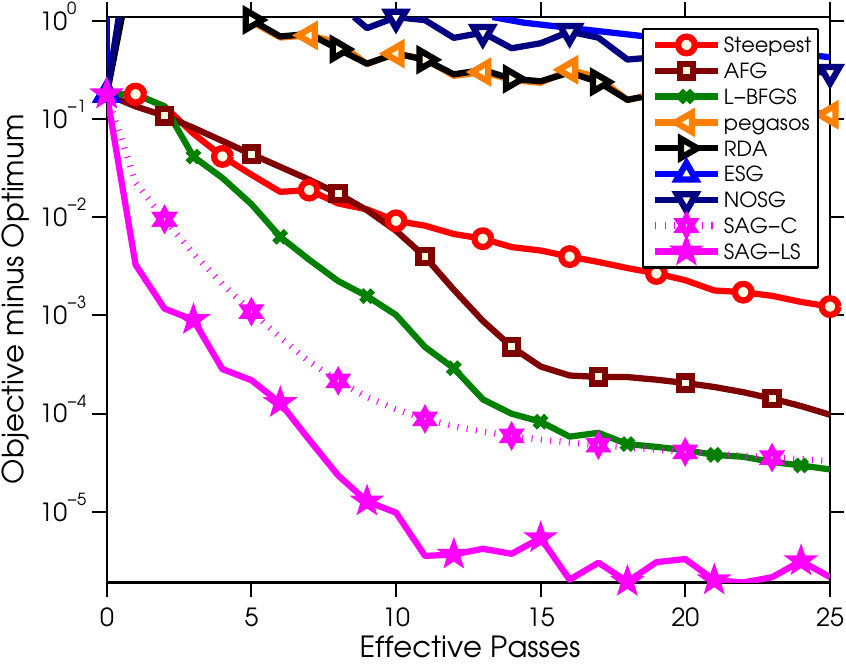} \hspace*{-.1cm}
\includegraphics[width=4.9cm]{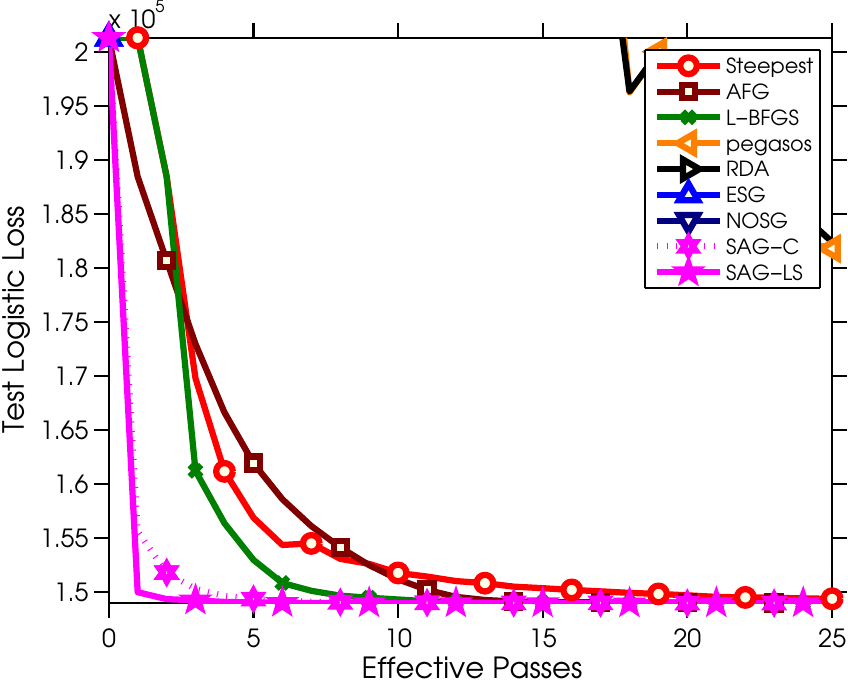} \hspace*{-.1cm}
\includegraphics[width=4.9cm]{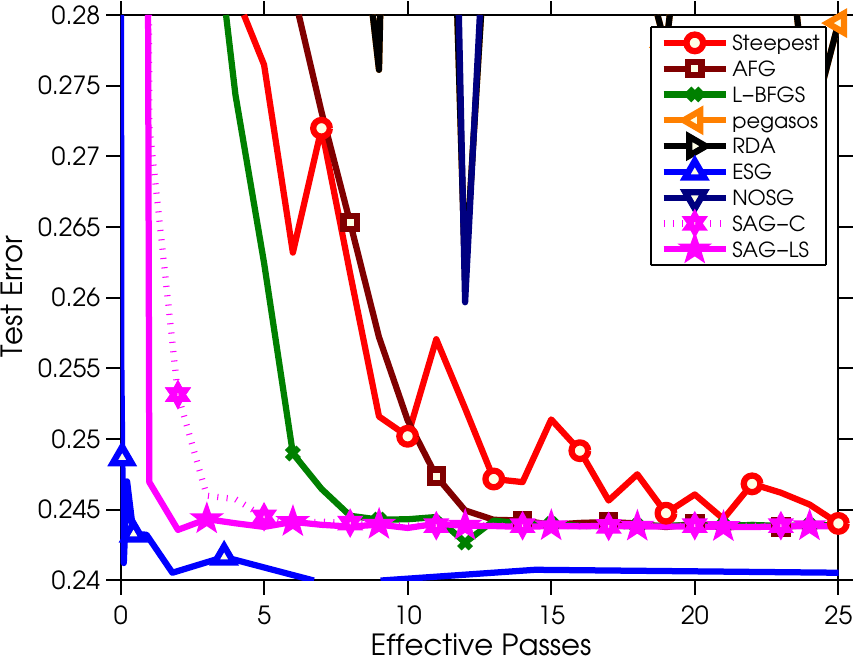}}

\vspace*{-.25cm}

\caption{Comparison of optimization strategies for $\ell_2$-regularized
logistic regression. Left: training excess objective. Middle: testing objective. Right: test errors. From top to bottom are the results on the 
 \emph{quantum}, \emph{protein}, \emph{sido}, \emph{rcv1},  and  \emph{covertype} data sets.  This figure is best viewed in colour.}
\label{fig:logreg}
\end{figure}

The theoretical convergence rates suggest the following strategies for
deciding on whether to use an FG or an SG method:
\begin{enumerate}
\item If we can only afford one pass through the data, then an SG method should
be used.
\item If we can afford to do many passes through the data (say, several
hundred), then an FG method should be used.
\end{enumerate}
We expect that the SAG iterations will be most useful between these two extremes, where we can afford to do more than one pass through the data but cannot afford to do enough passes to warrant using FG algorithms like L-BFGS.  To test whether this is indeed the case on real data sets, we performed experiments on a set of freely available benchmark binary classification data sets.  The 
\emph{quantum} ($p=50 000$, $p=78$) and 
\emph{protein}   ($n=145 751$, $p=74$) data set was obtained from the KDD Cup 2004 website,\footnote{\small \url{http://osmot.cs.cornell.edu/kddcup}} the \emph{sido} data set was obtained from the Causality Workbench website,\footnote{ \small\url{http://www.causality.inf.ethz.ch/home.php}}
 while the \emph{rcv1}  ($n= 20 242$, $p=47 236$) and \emph{covertype}   ($n=581 012 $, $p=54$) data sets were obtained from the LIBSVM data website.\footnote{ \small\url{http://www.csie.ntu.edu.tw/~cjlin/libsvmtools/datasets}}
Although our method can be applied to any differentiable function, on these data sets we focus on the $\ell_2$-regularized logistic regression problem, with $\lambda=1/n$.
We split each dataset in two, training on one half and testing on the other half.  
We added a (regularized) bias term to all data sets, and for dense features we standardized so that they would have a mean of zero and a variance of one.
We measure the training and testing costs as a function of the number of effective passes through the data, measured as the number of $f_i'$ evaluations divided by $n$. These results are thus independent of the practical implementation of the algorithms. We plot the training and testing costs of the different methods for 30 effective passes  through the data in Figure~\ref{fig:logreg}.  



\begin{figure}
\centering
\hspace*{-1cm}
\mbox{ \includegraphics[width=5.25cm]{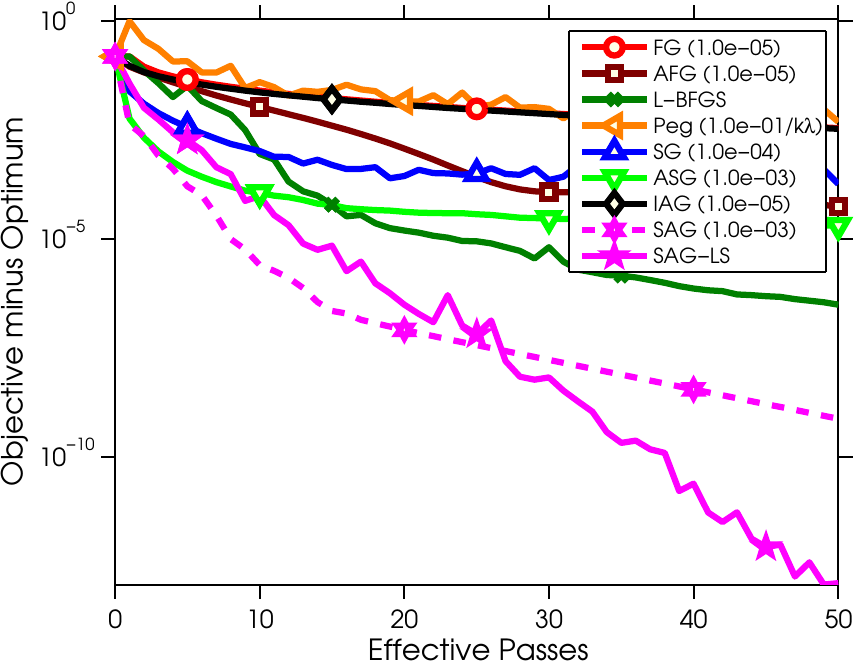} \hspace*{-.1cm}
\includegraphics[width=5.25cm]{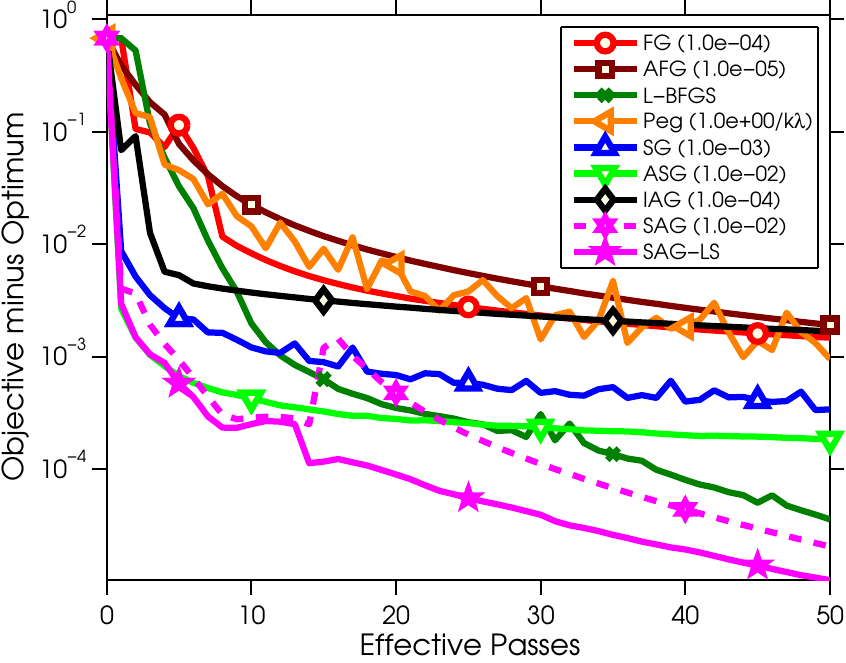} \hspace*{-.1cm}
\includegraphics[width=5.25cm]{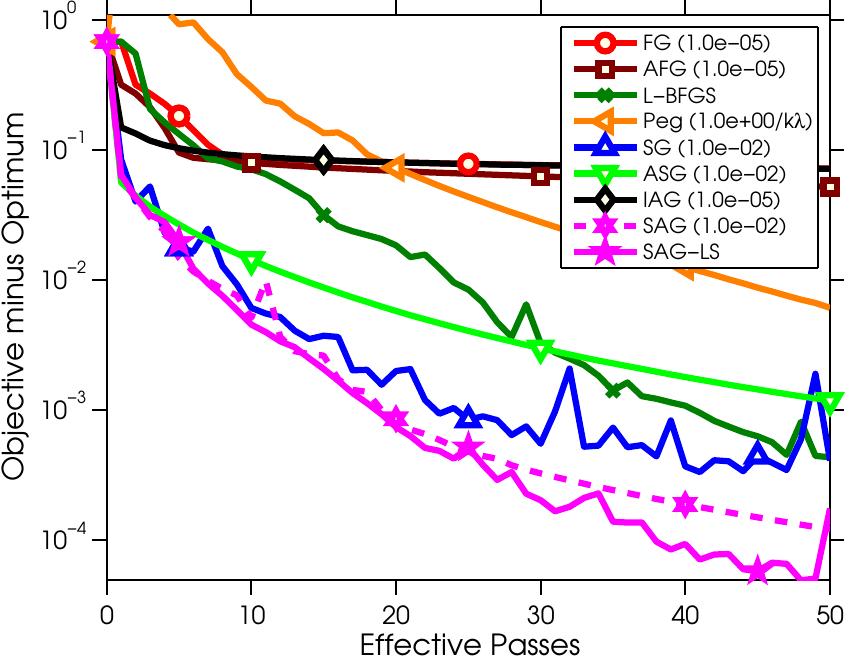} \hspace*{-.1cm}}

\mbox{\includegraphics[width=4.5cm]{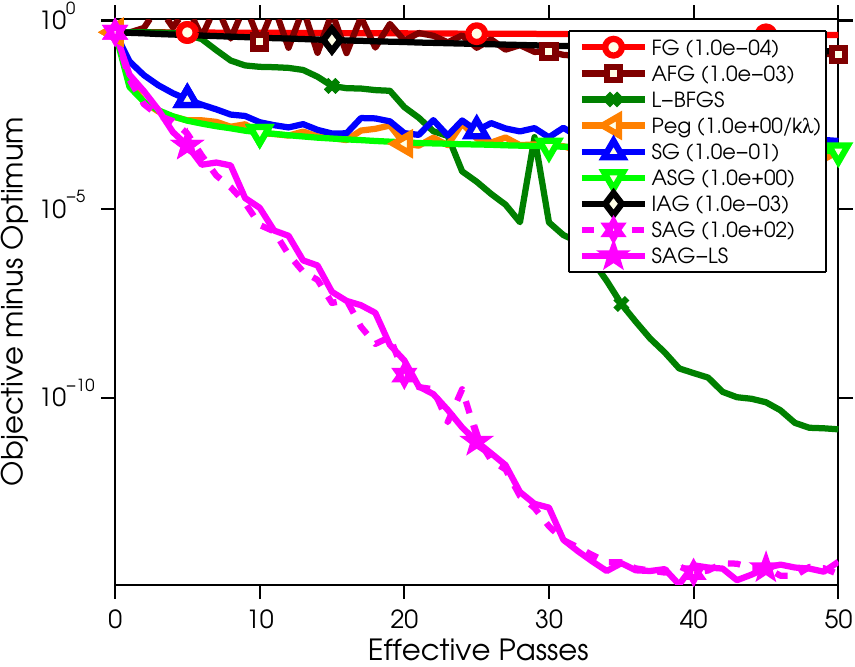} \hspace*{-.01cm}
\includegraphics[width=4.5cm]{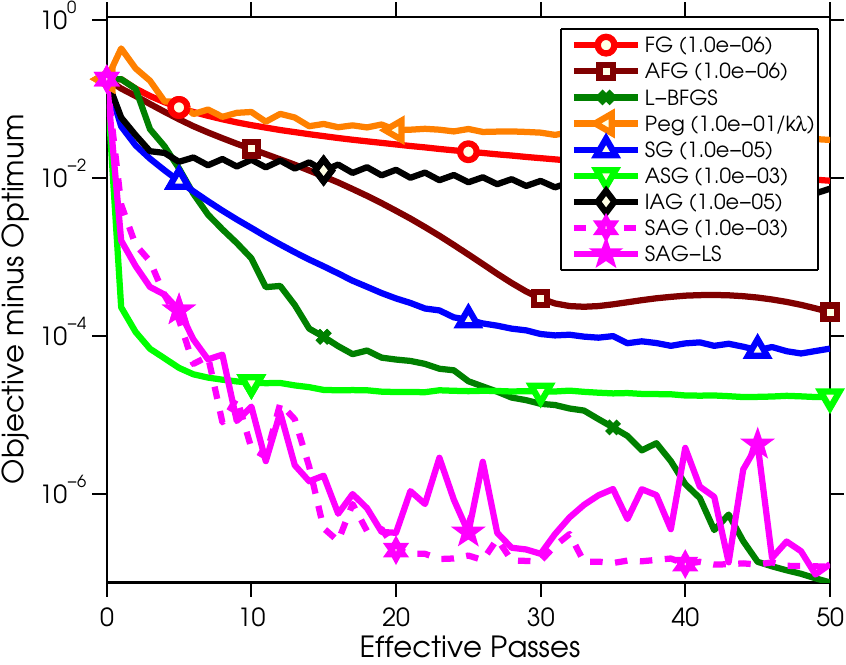} \hspace*{-.01cm}
}
\caption{Comparison of optimization strategies that choose the best step-size in hindsight. In the top row are the \emph{quantum} (left), \emph{protein} (center), and \emph{sido} (right) data sets. In the bottom row are the \emph{rcv1} and \emph{covertype} data sets. This figure is best viewed in colour.}
\label{fig:optim}
\end{figure}

In our second series of experiments, we sought to test whether SG methods (or the IAG method) with a very carefully chosen step size would be competitive with the SAG iterations.  In particular, we compared the following variety of basic FG and SG methods.
\begin{enumerate}
\item \textbf{FG}: The full gradient method described by
iteration~\eqref{eq:FG}.
\item \textbf{AFG}: The accelerated full gradient method~\cite{nesterov1983method}, where iterations of~\eqref{eq:FG} are
interleaved with an extrapolation step.
\item \textbf{peg}: The pegasos algorithm of~\cite{shalev2007pegasos}, but where we multiply the step size by a constant.
\item \textbf{SG}: The stochastic gradient method described by iteration~\eqref{eq:SG}, where we use a constant step-size.
\item \textbf{ASG}: The stochastic gradient method described by
iteration~\eqref{eq:SG}, where we use a constant step size and average the iterates.\footnote{We have also compared to a variety of other SG methods, such as SG with momentum, SG with gradient averaging, accelerated SG, and using SG but delaying averaging until after the first effective pass. However, none of these SG methods performed better than the ASG method above so we omit them to keep the plots simple.}
\item \textbf{IAG}: The incremental aggregated gradient method
of~\cite{blatt2008convergent} described by iteration~\eqref{eq:SAG} but with a
cyclic choice of $i_k$.
\item \textbf{SAG}: The proposed stochastic average gradient method
described by iteration~\eqref{eq:SAG}.
\end{enumerate}
For all of the above methods, we chose the step size that gave the best performance among powers of $10$.  On the full data sets, we compare these methods to each other and to the L-BFGS and the SAG-LS algorithms from the previous experiment in Figure~\ref{fig:optim}, which also shows the selected step sizes.

We can observe several trends across these experiments:
\begin{list}{\labelitemi}{\leftmargin=1.7em}
   \addtolength{\itemsep}{-.215\baselineskip}
\item[--] \textbf{FG vs. SG}: Although the performance of SG methods
 can be catastrophic if the step size is not chosen carefully (e.g., the \emph{quantum} and \emph{covertype} data), with a carefully-chosen step-size
the SG methods always do substantially better than FG methods
 on the first few passes through the data.  In contrast, the adaptive
 FG methods in the first experiment are not sensitive to the step size and because of its steady progress the best FG method (L-BFGS) 
 always eventually passes the SG methods.
\item[--] \textbf{(FG and SG) vs. SAG}: The SAG iterations seem to achieve the best of both worlds.  They
start out substantially better than FG methods, but continue
to make steady (linear) progress which leads to better performance than SG
methods.
The significant speed-up observed for SAG in reaching low training costs often also seems to translate into reaching the optimal testing cost more quickly than the other methods. We also note that the proposed line-search seems to perform as well or better than choosing the optimal fixed step-size in hind sight.

\item[--] \textbf{IAG vs. SAG}: The second experiment shows that the IAG method performs similarly to the regular FG method, and they also show the surprising result that the randomized SAG method outperforms the closely-related deterministic IAG method by a very large margin.  This is due to the larger step sizes used by the SAG iterations, which would cause the IAG iterations to diverge.
\end{list}




\section{Discussion}
\label{discussion}
\label{sec:discussion}

\textbf{Optimal regularization strength}: One might wonder if the additional hypothesis in Proposition~2 is satisfied in practice. In
a learning context, where each function $f_i$ is the loss
associated to a single data point, $L$ is equal to the largest
value of the loss second derivative $\xi$ (1 for the square loss, 1/4 for
the logistic loss) times $R^2$, where $R$ is a the uniform bound on the norm of each
data point. Thus, the constraint $\frac{\mu}{L}\geqslant \frac{8}{n}$ is
satisfied when
$\lambda \geqslant \frac{8 \xi R^2}{n}$. In low-dimensional settings, the optimal
regularization parameter is of the form $C/n$~\cite{liang} where
$C$ is a scalar constant, and may
thus violate the constraint. However, the improvement with respect to
regularization parameters of the form $\lambda =C/\sqrt{n}$ is known
to be asymptotically negligible, and in any case in such
low-dimensional settings, regular stochastic or batch
gradient descent may be efficient enough in practice.
In the more interesting high-dimensional settings where the
dimension $p$ of our covariates is not small compared to the sample
size $n$, then all theoretical analyses we are aware of advocate
settings of~$\lambda$ which satisfy this constraint.  For example,
\cite{srebro} considers parameters of the form $\lambda =
\frac{C}{\sqrt{n}}$ in the parametric setting, while \cite{steinwart}
considers  $\lambda = \frac{C}{n^{\beta}}$ with $\beta <1$ in a
non-parametric setting.

\textbf{Training cost vs.~testing cost}: The theoretical contribution of this work is limited to the convergence rate of the training cost. Though there are several settings where this is the metric of interest (e.g., variational inference in graphical models), in many cases one will be interested in the convergence speed of the testing cost. Since the $O(1/k)$ convergence rate of the testing cost, achieved by SG methods with decreasing step sizes (and a single pass through the data), is provably optimal when the algorithm only accesses the function through unbiased measurements of the objective and its gradient, it is unlikely that one can obtain a linear convergence rate for the testing cost with the SAG iterations. However, as shown in our experiments, the testing cost of the SAG iterates often reaches its minimum quicker than existing SG methods, and we could expect to improve the constant in the $O(1/k)$ convergence rate, as is the case with online second-order methods~\cite{bottou-bousquet-2011}.


 \textbf{Step-size selection and termination criteria}: The three major disadvantages of SG methods are: (i) the slow convergence rate, (ii) deciding when to terminate the algorithm, and (iii) choosing the step size while running the algorithm.  This paper showed that the SAG iterations achieve a much faster convergence rate, but the SAG iterations may also be advantageous in terms of tuning step sizes and designing termination criteria.  In particular, the SAG iterations suggest a natural termination criterion; 
 since the average of the $y_i^k$ variables converges to $g'(x^k)$ as $\norm{x^k-x^{k-1}}$ converges to zero, we can use $(1/n)\norm{\sum_iy_i^k}$ as an approximation of the optimality of $x^k$. Further, while SG methods require specifying a sequence of step sizes and mispecifying this sequence can have a disastrous effect on the convergence rate~\cite[\S2.1]{nemirovski2009robust}, our theory shows that the SAG iterations iterations achieve a linear convergence rate for any sufficiently small constant step size and our experiments indicate that a simple line-search gives strong performance.



\section*{Acknowledgements}


Nicolas Le Roux, Mark Schmidt, and Francis Bach are supported by the European Research Council (SIERRA-ERC-239993).  Mark Schmidt is also supported by a postdoctoral fellowship from the Natural Sciences and Engineering Research Council of Canada (NSERC).

\newpage

\appendix

\section*{Appendix}


In this Appendix, we first give the proofs of the two propositions. Subsequently, we compare the convergence rates of primal and dual FG and coordinate-wise methods to the rates of SAG for $\ell_2$-regularized least squares in terms of effective passes through the data.

\section{Proofs of the propositions}
We present here the proofs of Propositions~\ref{prop:small_alpha} and~\ref{prop:large_alpha}.

\subsection{Problem set-up and notations}
We use $g = \frac{1}{n}\sum_{i=1}^n f_i$ to denote a $\mu-$strongly convex function, where the functions $f_i$, $i=1,\ldots,n$ are convex functions from $\rb^p$ to $\rb$ with $L$-Lipschitz continuous gradients. Let us denote by $x^\ast$ the unique minimizer of $g$.

For $k\geqslant 1$, the stochastic average gradient algorithm performs the recursion
\BEAS
x^k & = & x^{k-1} - \frac{\alpha}{n} \sum_{i=1}^n y_i^{k},
\EEAS
where an $i_k$ is selected in $\{1,\dots,n\}$ uniformly at random and we set
\[
y^k_i = \begin{cases}
f_i'(x^{k-1}) & \textrm{if $i = i_k$,}\\
y^{k-1}_i & \textrm{otherwise.}
\end{cases}
\]

Denoting $z^k_i$ a random variable which takes the value $1-\frac{1}{n}$ with probability $\frac{1}{n}$ and $-\frac{1}{n}$ otherwise (thus with zero expectation), this is equivalent to
\BEAS
y_i^k & = & \left(1-\frac{1}{n}\right) y_i^{k-1} + \frac{1}{n}f'_{i}(x^{k-1}) + z^k_i \left[  f'_{i}(x^{k-1}) - y_i^{k-1}\right] \\
x^k & = & x^{k-1} - \frac{\alpha}{n} \sum_{i=1}^n \left[\left(1-\frac{1}{n}\right) y_i^{k-1} + \frac{1}{n}f'_{i}(x^{k-1}) + z^k_i \left[  f'_{i}(x^{k-1}) - y_i^{k-1}\right]\right]\\
&=& x^{k-1} - \frac{\alpha}{n}\left[\left(1-\frac{1}{n}\right)e^\top y^{k-1} + g'(x^{k-1}) + (z^k)^\top \left[ f'(x^{k-1}) - y^{k-1}\right]\right],
\EEAS
with
\begin{align*}
e = \left(
\begin{array}{c}
\idm\\
\vdots\\
\idm
\end{array}\right) \in \rb^{n p \times p},
\qquad
f'(x) = \left(
\begin{array}{c}
f_1'(x)\\
\vdots\\
f_n'(x)
\end{array}\right) \in \rb^{n p},
\qquad
z^k =
\left(
\begin{array}{c}
z_1^k I\\
\vdots\\
z_n^k I
\end{array}\right) \in \rb^{np \times p}.
\end{align*}
Using this definition of $z^k$, we have $\E [(z^k) (z^k)^\top] = \frac{1}{n} \idm - \frac{1}{n^2} ee^\top$. Note that, for a given $k$, the variables $z_1^k,\dots,z_n^k$ are not independent.

We also use the notation
\begin{align*}
\theta^k =
\left(
\begin{array}{c}
y_1^k\\
\vdots\\
y_n^k\\
x^k
\end{array}\right) \in \rb^{(n+1)p},
\qquad
\theta^\ast =
\left(
\begin{array}{c}
f_1'(x^\ast)\\
\vdots\\
f_n'(x^\ast)\\
x^\ast
\end{array}\right) \in \rb^{(n+1)p} \; .
\end{align*}

Finally, if $M$ is a $tp \times tp$ matrix and $m$ is a $tp \times p$ matrix, then:
\begin{itemize}
\item $\diag(M)$ is the $tp \times p$ matrix being the concatenation of the $t$ ($p \times p$)-blocks on the diagonal of $M$;
\item $\Diag(m)$ is the $tp \times tp$ block-diagonal matrix whose ($p \times p$)-blocks on the diagonal are equal to the ($p \times p$)-blocks of $m$.
\end{itemize}

\subsection{Outline of the proofs}
Each Proposition will be proved in multiple steps.
\begin{enumerate}
\item We shall find a Lyapunov function $Q$ from $\rb^{(n+1)p}$ to $\rb$ such that the sequence $\E Q(\theta^k)$ decreases at a linear rate.
\item We shall prove that $Q(\theta^k)$ dominates $\|x^k - x^\ast\|^2$ (in the case of Proposition~\ref{prop:large_alpha}) or $g(x^k) - g(x^\ast)$ (in the case of Proposition~\ref{prop:large_alpha}) by a constant for all $k$.
\item In the case of Proposition~\ref{prop:large_alpha}, we show how using one pass of stochastic gradient as the initialization provides the desired result.
\end{enumerate}

Throughout the proofs, $\mathcal{F}_{k}$ will denote the $\sigma$-field of information up to (and including time $k$), i.e., $\mathcal{F}_k$ is the $\sigma$-field generated by $z^1,\dots,z^k$.

\subsection{Convergence results for stochastic gradient descent}
\label{sec:bound_values_sg}
The constant in both our bounds depends on the initialization chosen. While this does not affect the linear convergence of the algorithm, the bound we obtain for the first few passes through the data is the $O(1/k)$ rate one would get using stochastic gradient descent, but with a constant proportional to $n$. This problem can be alleviated for the second bound by running stochastic gradient descent for a few iterations before running the SAG algorithm. In this section, we provide bounds for the stochastic gradient descent algorithm which will prove useful for the SAG algorithm.

The assumptions made in this section about the functions $f_i$ and the function $g$ are the same as the ones used for SAG. To get initial values for $x^0$ and $y^0$, we will do one pass of standard stochastic gradient.

We denote by $\sigma^2 = \frac{1}{n}\sum_{i=1}^n\|f_i'(x^\ast)\|^2$ the variance of the gradients at the optimum. We will use the following recursion:
$$
\tilde{x}^k = \tilde{x}^{k-1} - \gamma_k f_{i_k}'\left(\tilde{x}^{k-1}\right) \; .
$$

Denoting $\delta_k = \E \| \tilde{x}^k - x^\ast\|^2$, we have (following~\cite{bach2011non})
$$
\delta_k \leqslant \delta_{k-1} - 2 \gamma_k ( 1- \gamma_k L ) \E \left[  g'(\tilde{x}^{k-1})^\top( \tilde{x}^{k-1} - x^\ast) \right]+ 2 \gamma_k^2 \sigma^2 \; .
$$

Indeed, we have
\begin{align*}
\| \tilde{x}^k - x^\ast \|^2
& = \| \tilde{x}^{k-1} - x^\ast \|^2 - 2\gamma_k f_{i_k}'(\tilde{x}^{k-1})^\top ( \tilde{x}^{k-1} - x^\ast) + \gamma_k^2 \| f_{i_k}'(\tilde{x}^{k-1}) \|^2
\\
& \leqslant \| \tilde{x}^{k-1} - x^\ast \|^2 - 2\gamma_k f_{i_k}'(\tilde{x}^{k-1})^\top ( \tilde{x}^{k-1} - x^\ast) + 2 \gamma_k^2 \| f_{i_k}'(x^\ast) \|^2
+ 2 \gamma_k^2 \| f_{i_k}'(\tilde{x}^{k-1}) - f_{i_k}'(x^\ast) \|^2 \\
& \leqslant \| \tilde{x}^{k-1} - x^\ast \|^2 - 2\gamma_k f_{i_k}'(\tilde{x}^{k-1})^\top ( \tilde{x}^{k-1} - x^\ast) + 2\gamma_k^2 \| f_{i_k}'(x^\ast) \|^2\\
&\hspace*{.5cm} + 2 L\gamma_k^2 ( f_{i_k}'(\tilde{x}^{k-1}) - f_{i_k}'(x^\ast) ) ^\top ( \tilde{x}^{k-1} - x^\ast) \; .
\end{align*}
By taking expectations, we get
\BEAS
\E \left[ \| \tilde{x}^k - x^\ast \|^2 | \mathcal{F}_{k-1} \right]
& \leqslant &
\| \tilde{x}^{k-1} - x^\ast \|^2 - 2\gamma_k g'(\tilde{x}^{k-1})^\top ( \tilde{x}^{k-1} - x^\ast) + 2\gamma_k^2 \sigma^2
+ 2 L\gamma_k^2 g'(\tilde{x}^{k-1})^\top ( \tilde{x}^{k-1} - x^\ast) \\
\E \left[ \| \tilde{x}^k - x^\ast \|^2 \right]
& \leqslant &
\E \left[ \| \tilde{x}^{k-1} - x^\ast \|^2 \right] - 2\gamma_k ( 1 - \gamma_k L ) \E \left[ g'(\tilde{x}^{k-1})^\top ( \tilde{x}^{k-1} - x^\ast) \right] + 2\gamma_k^2 \sigma^2
\EEAS

Thus, if we take
$$
\gamma_k = \frac{1}{2L +  \frac{\mu}{2}k} \; ,
$$
we have $\gamma_k \leqslant 2\gamma_k ( 1 - \gamma_k L )$ and
\begin{align*}
\delta_k & \leqslant \delta_{k-1} -  \gamma_k \E \left[ g'(\tilde{x}^{k-1})^\top ( x^{k-1} - x^\ast) \right] + 2 \gamma_k^2 \sigma^2\\
& \leqslant \delta_{k-1} -  \gamma_k \left[ \E \left[ g(x^{k-1}) - g(x^\ast)\right] + \frac{\mu}{2} \delta_{k-1} \right]   + 2 \gamma_k^2 \sigma^2 \textrm{ using the strong convexity of $g$}\\
\E g(x^{k-1}) - g(x^\ast) & \leqslant - \frac{1}{ \gamma_k} \delta_{k} + \left( \frac{1}{\gamma_k} - \frac{\mu}{2}\right) \delta_{k-1} + 2 \gamma_k \sigma^2 \\
& \leqslant - \left( 2L + \frac{\mu}{2}k\right) \delta_{k} + \left( 2L + \frac{\mu}{2}(k-1)\right) \delta_{k-1} + 2 \gamma_k \sigma^2 \; .
\end{align*}

Averaging from $i = 0$ to $k-1$ and using the convexity of $g$, we have
\begin{align*}
\frac{1}{k}\sum_{i=0}^{k-1}\E g(x^{k-1}) - g(x^\ast) & \leqslant \frac{2L}{k}\delta_{0} + \frac{2 \sigma^2 }{k} \sum_{i=1}^k \gamma_i \\
\E g \left( \frac{1}{k}\sum_{i=0}^{k-1} x^i \right) - g(x^\ast)
& \leqslant \frac{2L}{k} \delta_0 + \frac{2 \sigma^2 }{k} \sum_{i=1}^k \gamma_i \\
& \leqslant \frac{2L}{k}\| x^0 - x^\ast\|^2 + \frac{2 \sigma^2 }{k} \sum_{i=1}^k \frac{1}{2L +  \frac{\mu}{2}i}\\
& \leqslant \frac{2L}{k}L \| x^0 - x^\ast\|^2 + \frac{2 \sigma^2 }{k} \int_0^k \frac{1}{2L +  \frac{\mu}{2}t}dt\\
& \leqslant \frac{2L}{k}\| x^0 - x^\ast\|^2 + \frac{4 \sigma^2 }{k \mu} \log\left( 1 + \frac{\mu k }{4L } \right) \; .
\end{align*}

\subsection{Important lemma}

In both proofs, our Lyapunov function contains a quadratic term \mbox{$R(\theta^k) = (\theta^k - \theta^*)^\top \left(\begin{array}{cc} A & b \\ b^\top & c \end{array}\right) (\theta^k - \theta^*)$} for some values of $A$, $b$ and $c$. The lemma below computes the value of $R(\theta^k)$ in terms of elements of $\theta^{k-1}$.

\begin{lemma}
If $ P = \left(
\begin{array}{cc} A & b \\ b^\top & c
\end{array}
\right)$,
for $A \in \rb^{np \times np}$, $b \in \rb^{np \times p}$ and $c \in \rb^{p \times p}$, then
\begin{align*}
&\E\left[\left.(\theta^k - \theta^*)^\top \left(\begin{array}{cc} A & b \\ b^\top & c
\end{array}\right) (\theta^k - \theta^*)\right| \mathcal{F}_{k-1}\right] \nonumber\\
&\hspace*{.5cm}=(y^{k-1} - f'(x^\ast))^\top \left[\left(1 - \frac{2}{n}\right)S + \frac{1}{n}\Diag(\diag(S))\right](y^{k-1} - f'(x^\ast))\\
&\hspace*{.5cm}+ \frac{1}{n}(f'(x^{k-1}) - f'(x^\ast))^\top \Diag(\diag(S))(f'(x^{k-1}) - f'(x^\ast))\\
&\hspace*{.5cm}+ \frac{2}{n}(y^{k-1} - f'(x^\ast))^\top \left[S - \Diag(\diag(S))\right] (f'(x^{k-1}) - f'(x^\ast))\\
&\hspace*{.5cm}+ 2\left(1 - \frac{1}{n}\right) (y^{k-1} - f'(x^\ast))^\top \left[b - \frac{\alpha}{n}ec\right](x^{k-1} - x^\ast)\\
&\hspace*{.5cm}+ \frac{2}{n}(f'(x^{k-1}) - f'(x^\ast))^\top\left[b - \frac{\alpha}{n}ec\right](x^{k-1} - x^\ast)\\
&\hspace*{.5cm}+ (x^{k-1} - x^\ast)^\top c (x^{k-1} - x^\ast) \; ,
\end{align*}
with
\[
S = A - \frac{\alpha}{n} b e^\top - \frac{\alpha}{n} e b^\top + \frac{\alpha^2}{n^2} e ce^\top \; .
\]
\end{lemma}
Note that for square $n \times n$ matrix, $\diag(M)$ denotes a vector of size $n$ composed of the diagonal of $M$, while for a vector $m$ of dimension $n$,
$\Diag(m)$ is the $n \times n$ diagonal matrix with $m$ on its diagonal. Thus $\Diag(\diag(M))$ is a diagonal matrix with the diagonal elements of $M$ on its diagonal, and $\diag(\Diag(m))=m$.

\begin{proof}
Throughout the proof, we will use the equality $g'(x) = e^\top f'(x) / n$. Moreover, all conditional expectations of linear functions of $z^k$ will be equal to zero.

We have
\begin{align}
&\E\left[\left.(\theta^k - \theta^*)^\top \left(\begin{array}{cc} A & b \\ b^\top & c
\end{array}\right) (\theta^k - \theta^*)\right| \mathcal{F}_{k-1}\right] \nonumber\\
&\hspace*{.5cm}= E\left[ (y^k - f'(x^\ast))^\top A (y^k - f'(x^\ast))+ 2 (y^k - f'(x^\ast))^\top b (x^k - x^\ast) + (x^k - x^\ast)^\top c (x^k - x^\ast)| \mathcal{F}_{k-1}\right] \; .
\label{eq:recursion}
\end{align}

The first term (within the expectation) on the right-hand side of \eq{recursion} is equal to
\begin{align*}
(y^k - f'(x^\ast))^\top A (y^k - f'(x^\ast)) & = \left(1 - \frac{1}{n}\right)^2(y^{k-1} - f'(x^\ast))^\top A (y^{k-1} - f'(x^\ast))\\
& \hspace*{.5cm} + \frac{1}{n^2}(f'(x^{k-1}) - f'(x^\ast))^\top A (f'(x^{k-1}) - f'(x^\ast))\\
& \hspace*{.5cm} + [\Diag(z^k) (f'(x^{k-1}) - y^{k-1})]^\top A [\Diag(z^k) (f'(x^{k-1}) - y^{k-1})]\\
& \hspace*{.5cm} + \frac{2}{n}\left(1 - \frac{1}{n}\right)(y^{k-1} - f'(x^\ast))^\top A (f'(x^{k-1}) - f'(x^\ast)) \; .
\end{align*}

The only random term (given $\mathcal{F}_{k-1}$) is the third one whose expectation is equal to
\begin{align*}
&\E\left[[\Diag(z^k) (f'(x^{k-1}) - y^{k-1})]^\top A [\Diag(z^k) (f'(x^{k-1}) - y^{k-1})]| \mathcal{F}_{k-1}\right]\\
&\hspace*{.5cm} = \frac{1}{n}(f'(x^{k-1}) - y^{k-1})^\top\left[\Diag(\diag(A)) - \frac{1}{n}A\right](f'(x^{k-1}) - y^{k-1}) \; .
\end{align*}

The second term (within the expectation) on the right-hand side of \eq{recursion} is equal to
\begin{align*}
(y^k - f'(x^\ast))^\top b (x^k - x^\ast) 	& = \left(1 - \frac{1}{n}\right)(y^{k-1} - f'(x^\ast))^\top b (x^{k-1} - x^\ast) \\
& \hspace*{.5cm} + \frac{1}{n}(f'(x^{k-1}) - f'(x^\ast))^\top b (x^{k-1} - x^\ast) \\
& \hspace*{.5cm} - \frac{\alpha}{n} \left(1 - \frac{1}{n}\right)^2(y^{k-1} - f'(x^\ast))^\top b e^\top (y^{k-1} - f'(x^\ast))\\
& \hspace*{.5cm} - \frac{\alpha}{n} \frac{1}{n} \left(1 - \frac{1}{n}\right)(f'(x^{k-1}) - f'(x^\ast))^\top b e^\top (y^{k-1} - f'(x^\ast))\\
& \hspace*{.5cm} - \frac{\alpha}{n} \frac{1}{n} \left(1 - \frac{1}{n}\right)(y^{k-1} - f'(x^\ast))^\top b e^\top (f'(x^{k-1}) - f'(x^\ast))\\
& \hspace*{.5cm} - \frac{\alpha}{n} \frac{1}{n^2} (f'(x^{k-1}) - f'(x^\ast))^\top b e^\top(f'(x^{k-1}) - f'(x^\ast))\\
& \hspace*{.5cm} - \frac{\alpha}{n} [\Diag(z^k) (f'(x^{k-1}) - y^{k-1})]^\top b (z^k)^\top \left[(f'(x^{k-1}) - y^{k-1})\right]
\end{align*}

The only random term (given $\mathcal{F}_{k-1}$) is the last one whose expectation is equal to
\begin{align*}
&\E\left[[\Diag(z^k) (f'(x^{k-1}) - y^{k-1})]^\top b (z^k)^\top \left[(f'(x^{k-1}) - y^{k-1})\right]| \mathcal{F}_{k-1}\right]\\
&\hspace*{.5cm} = \frac{1}{n}(f'(x^{k-1}) - y^{k-1})^\top \left(\Diag(\diag(be^\top) - \frac{1}{n}b e^\top\right)(f'(x^{k-1}) - y^{k-1}) \; .
\end{align*}

The last term on the right-hand side of \eq{recursion} is equal to
\begin{align*}
(x^k - x^\ast)^\top c (x^k - x^\ast) &= (x^{k-1} - x^\ast)^\top c (x^{k-1} - x^\ast)\\
& \hspace*{.5cm} + \frac{\alpha^2}{n^2}\left(1 - \frac{1}{n}\right)^2 (y^{k-1} - f'(x^\ast))^\top e c e^\top (y^{k-1} - f'(x^\ast))\\
& \hspace*{.5cm} + \frac{\alpha^2}{n^2} \frac{1}{n^2} (f'(x^{k-1}) - f'(x^\ast))^\top e c e^\top (f'(x^{k-1}) - f'(x^\ast))\\
& \hspace*{.5cm} -\frac{2\alpha}{n}\left(1 - \frac{1}{n}\right) (x^{k-1} - x^\ast)^\top c e^\top (y^{k-1} - f'(x^\ast))\\
& \hspace*{.5cm} -\frac{2\alpha}{n}\frac{1}{n}(x^{k-1} - x^\ast)^\top c e^\top (f'(x^{k-1}) - f'(x^\ast))\\
& \hspace*{.5cm} + \frac{2\alpha^2}{n^2} \frac{1}{n}\left(1 - \frac{1}{n}\right)(y^{k-1} - f'(x^\ast))^\top e c e^\top (f'(x^{k-1}) - f'(x^\ast))\\
& \hspace*{.5cm} + \frac{\alpha^2}{n^2} \left[(z^k)^\top(f'(x^{k-1}) - y^{k-1})\right]^\top c \left[(z^k)^\top(f'(x^{k-1}) - y^{k-1})\right] \; .
\end{align*}
The only random term (given $\mathcal{F}_{k-1}$) is the last one whose expectation is equal to
\begin{align*}
&\E\left[ \left[(z^k)^\top(f'(x^{k-1}) - y^{k-1})\right]^\top c \left[(z^k)^\top(f'(x^{k-1}) - y^{k-1})\right]| \mathcal{F}_{k-1}\right]\\
&\hspace*{.5cm} = \frac{1}{n}(f'(x^{k-1}) - y^{k-1})^\top \left[\Diag(\diag(ece^\top)) - \frac{1}{n}ece^\top\right](f'(x^{k-1}) - y^{k-1}) \; .
\end{align*}

Summing all these terms together, we get the following result:
\begin{align*}
&\E\left[\left.(\theta^k - \theta^*)^\top \left(\begin{array}{cc} A & b \\ b^\top & c
\end{array}\right) (\theta^k - \theta^*)\right| \mathcal{F}_{k-1}\right] \nonumber\\
&\hspace*{.5cm}=\left(1 - \frac{1}{n}\right)^2 (y^{k-1} - f'(x^\ast))^\top S(y^{k-1} - f'(x^\ast))\\
&\hspace*{.5cm}+ \frac{1}{n^2}(f'(x^{k-1}) - f'(x^\ast))^\top S(f'(x^{k-1}) - f'(x^\ast))\\
&\hspace*{.5cm}+ \frac{1}{n}(f'(x^{k-1}) - y^{k-1})^\top\left[\Diag(\diag(S)) - \frac{1}{n}S\right](f'(x^{k-1}) - y^{k-1})\\
&\hspace*{.5cm}+ \frac{2}{n}\left(1 - \frac{1}{n}\right)(y^{k-1} - f'(x^\ast))^\top S (f'(x^{k-1}) - f'(x^\ast))\\
&\hspace*{.5cm}+ 2\left(1 - \frac{1}{n}\right) (y^{k-1} - f'(x^\ast))^\top \left[b - \frac{\alpha}{n}ec\right](x^{k-1} - x^\ast)\\
&\hspace*{.5cm}+ \frac{2}{n}(f'(x^{k-1}) - f'(x^\ast))^\top\left[b - \frac{\alpha}{n}ec\right](x^{k-1} - x^\ast)\\
&\hspace*{.5cm}+ (x^{k-1} - x^\ast)^\top c (x^{k-1} - x^\ast)
\end{align*}
with $S = A - \frac{\alpha}{n} b e^\top - \frac{\alpha}{n} e b^\top + \frac{\alpha^2}{n^2} e ce^\top
= A - b c^{-1} b^\top + (b - \frac{\alpha}{n} e c ) c^{-1} (b - \frac{\alpha}{n} e c )^\top
$.

Rewriting $f'(x^{k-1}) - y^{k-1} = (f'(x^{k-1}) - f'(x^\ast)) - (y^{k-1} - f'(x^\ast))$, we have
\begin{align*}
&f'(x^{k-1}) - y^{k-1})^\top\left[\Diag(\diag(S)) - \frac{1}{n}S\right](f'(x^{k-1}) - y^{k-1})\\
&\hspace*{.5cm} = (f'(x^{k-1}) - f'(x^\ast))^\top\left[\Diag(\diag(S)) - \frac{1}{n}S\right](f'(x^{k-1}) - f'(x^\ast))\\
&\hspace*{.5cm} + (y^{k-1} - f'(x^\ast))^\top\left[\Diag(\diag(S)) - \frac{1}{n}S\right](y^{k-1} - f'(x^\ast))\\
&\hspace*{.5cm} - 2(y^{k-1} - f'(x^\ast))^\top\left[\Diag(\diag(S)) - \frac{1}{n}S\right](f'(x^{k-1}) - f'(x^\ast)).
\end{align*}

Hence, the sum may be rewritten as
\begin{align*}
&\E\left[\left.(\theta^k - \theta^*)^\top \left(\begin{array}{cc} A & b \\ b^\top & c
\end{array}\right) (\theta^k - \theta^*)\right| \mathcal{F}_{k-1}\right] \nonumber\\
&\hspace*{.5cm}=(y^{k-1} - f'(x^\ast))^\top \left[\left(1 - \frac{2}{n}\right)S + \frac{1}{n}\Diag(\diag(S))\right](y^{k-1} - f'(x^\ast))\\
&\hspace*{.5cm}+ \frac{1}{n}(f'(x^{k-1}) - f'(x^\ast))^\top \Diag(\diag(S))(f'(x^{k-1}) - f'(x^\ast))\\
&\hspace*{.5cm}+ \frac{2}{n}(y^{k-1} - f'(x^\ast))^\top \left[S - \Diag(\diag(S))\right] (f'(x^{k-1}) - f'(x^\ast))\\
&\hspace*{.5cm}+ 2\left(1 - \frac{1}{n}\right) (y^{k-1} - f'(x^\ast))^\top \left[b - \frac{\alpha}{n}ec\right](x^{k-1} - x^\ast)\\
&\hspace*{.5cm}+ \frac{2}{n}(f'(x^{k-1}) - f'(x^\ast))^\top\left[b - \frac{\alpha}{n}ec\right](x^{k-1} - x^\ast)\\
&\hspace*{.5cm}+ (x^{k-1} - x^\ast)^\top c (x^{k-1} - x^\ast)
\end{align*}

This concludes the proof.
\end{proof}

\subsection{Analysis for $\alpha = \frac{1}{2nL}$ }
We now prove Proposition~\ref{prop:small_alpha}, providing a bound for the convergence rate of the SAG algorithm in the case of a small step size, $\alpha = \frac{1}{2nL}$.

\begin{proof}
\subsubsection*{Step 1 - Linear convergence of the Lyapunov function}
In this case, our Lyapunov function is quadratic, i.e.,
\[
Q(\theta^k) = (\theta^k - \theta^*)^\top \left(\begin{array}{cc} A & b \\ b^\top & c \end{array}\right) (\theta^k - \theta^*) \; .
\]

We consider
\BEAS
A & = & 3n \alpha^2 \idm + \frac{\alpha^2}{n} (\frac{1}{n} - 2 )ee^\top \\
b & = & - \alpha ( 1 - \frac{1}{n}) e \\
c & = & \idm \\
S & = & 3 n \alpha^2 \idm \\
b - \frac{\alpha}{n} e c & = & - \alpha e \; .
\EEAS

The goal will be to prove that $\E [ Q(\theta^k) | \mathcal{F}_{k-1} ] - ( 1 - \delta)Q(\theta^{k-1})$ is negative for some $\delta > 0$. This will be achieved by bounding all the terms by a term depending on $g'(x^{k-1})^\top( x^{k-1} - x^{\ast})$ whose positivity is guaranteed by the convexity of $g$.

We have, with our definition of $A$, $b$ and $c$:
\begin{align*}
S - \Diag(\diag(S)) &= 3 n \alpha^2 \idm - 3 n \alpha^2 \idm = 0\\
e^\top (f'(x^{k-1}) - f'(x^\ast)) &= n[g'(x^{k-1}) - g'(x^\ast)] = ng'(x^{k-1}) \; .
\end{align*}

This leads to (using the lemma of the previous section):
\begin{align*}
\E [ Q(\theta^k) | \mathcal{F}_{k-1} ] &=
\E \bigg[ (\theta^k - \theta^\ast)^\top
\left(
\begin{array}{cc} A & b \\ b^\top & c
\end{array}
\right)
(\theta^k - \theta^\ast) \bigg| \mathcal{F}_{k-1}\bigg]
\\
& =
\left(1-\frac{1}{n}\right) 3 n\alpha^2 (y^{k-1} - f'(x^\ast))^\top
(y^{k-1} - f'(x^\ast))
\\
&\hspace*{.5cm} + (x^{k-1} -x^\ast)^\top
(x^{k-1} -x^\ast) - \frac{2 \alpha}{n} (x^{k-1} -x^\ast)^\top e^\top ( f'(x^{k-1}) - f'(x^\ast)) \\
& \hspace*{.5cm}+ 3 \alpha^2 ( f'(x^{k-1}) - f'(x^\ast))^\top ( f'(x^{k-1}) - f'(x^\ast)) \\
&\hspace*{.5cm} - 2\alpha \left(1-\frac{1}{n}\right)(y^{k-1} - f'(x^\ast))^\top  e (x^{k-1} -x^\ast) \\
\\
& =
\left(1-\frac{1}{n}\right) 3 n\alpha^2 (y^{k-1} - f'(x^\ast))^\top
(y^{k-1} - f'(x^\ast))
\\
&\hspace*{.5cm} + (x^{k-1} -x^\ast)^\top
(x^{k-1} -x^\ast) - {2 \alpha} (x^{k-1} -x^\ast)^\top g'(x^{k-1}) \\
& \hspace*{.5cm}+ 3 \alpha^2 ( f'(x^{k-1}) - f'(x^\ast))^\top ( f'(x^{k-1}) - f'(x^\ast)) \\
&\hspace*{.5cm} - 2\alpha \left(1-\frac{1}{n}\right)(y^{k-1} - f'(x^\ast))^\top  e (x^{k-1} -x^\ast) \\
\\
&\leqslant
\left(1-\frac{1}{n}\right) 3 n\alpha^2 (y^{k-1} - f'(x^\ast))^\top
(y^{k-1} - f'(x^\ast))
\\
& \hspace*{.5cm}+ (x^{k-1} -x^\ast)^\top
(x^{k-1} -x^\ast) - {2 \alpha} (x^{k-1} -x^\ast)^\top g'(x^{k-1}) \\
&\hspace*{.5cm}+ 3 \alpha^2 n L (x^{k-1} -x^\ast)^\top g'(x^{k-1}) \\
&\hspace*{.5cm} - 2\alpha \left(1-\frac{1}{n}\right)(y^{k-1} - f'(x^\ast))^\top  e (x^{k-1} -x^\ast) \; .
\end{align*}
The third line is obtained using the Lipschitz property of the gradient, that is
\begin{align*}
( f'(x^{k-1}) - f'(x^\ast))^\top ( f'(x^{k-1}) - f'(x^\ast))
& = \sum_{i=1}^n \| f_i'(x^{k-1}) - f_i'(x^{\ast})\|^2 \\
&
\leqslant  \sum_{i=1}^n L ( f_i'(x^{k-1}) - f_i'(x^{\ast}))^\top( x^{k-1} - x^{\ast})
\\
& = nL ( g'(x^{k-1}) - g'(x^{\ast}))^\top( x^{k-1} - x^{\ast}) \; ,
\end{align*}
where the inequality in the second line stems from~\cite[Theorem 2.1.5]{nesterov2004introductory}.

We have
\begin{align*}
(1 - \delta)Q(\theta^{k-1}) &= (1 - \delta)(\theta^{k-1} - \theta^\ast)^\top
\left(
\begin{array}{cc} A & b \\ b^\top & c
\end{array}
\right)
(\theta^{k-1} - \theta^\ast)\\
&= ( 1- \delta)  (y^{k-1} - f'(x^\ast))^\top
\left[ 3n \alpha^2 \idm + \frac{\alpha^2}{n} \left(\frac{1}{n} - 2 \right)ee^\top \right]
(y^{k-1} - f'(x^\ast)) \\
& \hspace*{.5cm} + ( 1 - \delta) (x^{k-1} -x^\ast)^\top
(x^{k-1} -x^\ast)\\
& \hspace*{.5cm} - 2\alpha ( 1 - \delta) \left(1-\frac{1}{n}\right)(y^{k-1} - f'(x^\ast))^\top  e (x^{k-1} -x^\ast) \; .
\end{align*}

The difference is then:
\begin{align*}
&\E [ Q(\theta^k) | \mathcal{F}_{k-1} ]  - (1 - \delta)Q(\theta^{k-1})\\
&\hspace*{.5cm}\leqslant (y^{k-1} - f'(x^\ast))^\top \left[
3 n \alpha^2 \left( \delta - \frac{1}{n} \right) \idm
+ ( 1 - \delta) \frac{\alpha^2}{n} \left( 2 - \frac{1}{n} \right) ee^\top
\right] (y^{k-1} - f'(x^\ast))
\\
& \hspace*{1cm}+ \delta (x^{k-1} -x^\ast)^\top
(x^{k-1} -x^\ast)\\
&\hspace*{1cm} - ( {2 \alpha} - 3 \alpha^2 n L ) (x^{k-1} -x^\ast)^\top g'(x^{k-1}) \\
&\hspace*{1cm} - 2\alpha \delta \left(1-\frac{1}{n}\right)(y^{k-1} - f'(x^\ast))^\top  e (x^{k-1} -x^\ast).
\end{align*}

Note that for any symmetric negative definite matrix $M$ and for any vectors $s$ and $t$ we have
\[
(s + \frac{1}{2}M^{-1}t)^\top M (s + \frac{1}{2}M^{-1}t) \leqslant 0,
\]
and thus that
\[
s^\top M s + s^\top t \leqslant -\frac{1}{4} t^\top M^{-1} t \; .
\]
Using this fact with
\begin{align*}
M &= \left[
3 n \alpha^2 \left( \delta - \frac{1}{n} \right) \idm
+ ( 1 - \delta) \frac{\alpha^2}{n} \left( 2 - \frac{1}{n} \right) ee^\top
\right]\\
&=
\left[
3 n \alpha^2 \left( \delta - \frac{1}{n} \right) \left( \idm - \frac{ee^\top}{n}\right)
+ \alpha^2\left(3n\delta - 1- 2\delta + \frac{\delta - 1}{n}\right) \frac{ee^\top}{n}
\right]\\
s &= y^{k-1} - f'(x^\ast)\\
t &= - 2\alpha \delta \left(1-\frac{1}{n}\right)e (x^{k-1} -x^\ast) \; ,
\end{align*}
we have
\begin{align*}
&(y^{k-1} - f'(x^\ast))^\top \left[
3 n \alpha^2 \left( \delta - \frac{1}{n} \right) \idm
+ ( 1 - \delta) \frac{\alpha^2}{n} \left( 2 - \frac{1}{n} \right) ee^\top
\right] (y^{k-1} - f'(x^\ast))\\
& \hspace*{1cm} - 2\alpha \delta \left(1-\frac{1}{n}\right)(y^{k-1} - f'(x^\ast))^\top  e (x^{k-1} -x^\ast)\\
& \hspace*{.5cm} \leqslant -\alpha^2\delta^2\left(1-\frac{1}{n}\right)^2 (x^{k-1} -x^\ast)^\top e^\top \left[
3 n \alpha^2 \left( \delta - \frac{1}{n} \right) \left( \idm - \frac{ee^\top}{n}\right) \right. \\
& \hspace*{1.5cm}
\left.+ \alpha^2\left(3n\delta - 1- 2\delta + \frac{\delta - 1}{n}\right) \frac{ee^\top}{n}
\right]^{-1} e (x^{k-1} -x^\ast)\\
&\hspace*{.5cm}= -\frac{\alpha^2\delta^2\left(1-\frac{1}{n}\right)^2n}{\alpha^2\left[3n\delta - 1- 2\delta + \frac{\delta - 1}{n}\right]} \|x^{k-1} -x^\ast\|^2\\
&\hspace*{.5cm}= -\frac{\delta^2\left(1-\frac{1}{n}\right)^2n}{3n\delta - 1- 2\delta + \frac{\delta - 1}{n}} \|x^{k-1} -x^\ast\|^2 \; .
\end{align*}

A 
sufficient condition for $M$ to be negative definite is to have $\delta \leqslant \frac{1}{3n}$.

The bound then becomes
\begin{align*}
\E [ Q(\theta^k) | \mathcal{F}_{k-1} ]  - (1 - \delta)Q(\theta^{k-1}) &\leqslant - ( 2 \alpha - 3 \alpha^2 n L ) (x^{k-1} -x^\ast)^\top g'(x^{k-1})\\
&\hspace*{.5cm} + \left(\delta - \frac{\delta^2\left(1-\frac{1}{n}\right)^2}{\left[3n\delta - 1- 2\delta + \frac{\delta - 1}{n}\right]}n\right) \|x^{k-1} -x^\ast\|^2 \; .
\end{align*}

We now use the strong convexity of $g$ to get the inequality
\[
\|x^{k-1} -x^\ast\|^2 \leqslant \frac{1}{\mu}(x^{k-1} -x^\ast)^\top g'(x^{k-1}) \; .
\]

This yields the final bound
\begin{align*}
\E [ Q(\theta^k) | \mathcal{F}_{k-1} ]  - (1 - \delta)Q(\theta^{k-1}) &\leqslant - \left( 2 \alpha - 3 \alpha^2 n L + \frac{\delta^2\left(1-\frac{1}{n}\right)^2}{\left[3n\delta - 1- 2\delta + \frac{\delta - 1}{n}\right]}\frac{n}{\mu} - \frac{\delta}{\mu} \right) (x^{k-1} -x^\ast)^\top g'(x^{k-1}).
\end{align*}

Since we know that $(x^{k-1} -x^\ast)^\top g'(x^{k-1})$ is positive, due to the convexity of $g$, we need to prove that $\displaystyle \left( 2 \alpha - 3 \alpha^2 n L + \frac{\delta^2\left(1-\frac{1}{n}\right)^2}{\left[3n\delta - 1- 2\delta + \frac{\delta - 1}{n}\right]}\frac{n}{\mu} - \frac{\delta}{\mu} \right)$ is positive.

Using $\delta = \frac{\mu}{8 n L}$ and $\alpha = \frac{1}{2nL}$ gives
\begin{align*}
2 \alpha - 3 \alpha^2 n L + \frac{\delta^2\left(1-\frac{1}{n}\right)^2}{\left[3n\delta - 1- 2\delta + \frac{\delta - 1}{n}\right]}\frac{n}{\mu} - \frac{\delta}{\mu} &= \frac{1}{nL} - \frac{3}{4nL} - \frac{1}{8nL} - \frac{\delta^2\left(1-\frac{1}{n}\right)^2\frac{n}{\mu}}{1 - 3n\delta + 2\delta + \frac{1 - \delta}{n}}\\
&\geqslant \frac{1}{8nL} - \frac{\delta^2\frac{n}{\mu}}{1 - 3n\delta}\\
&= \frac{1}{8nL} - \frac{\frac{\mu}{64nL^2}}{1 - \frac{3\mu}{8L}}\\
&\geqslant \frac{1}{8nL} - \frac{\frac{\mu}{64nL^2}}{1 - \frac{3}{8}}\\
&= \frac{1}{8nL} - \frac{\mu}{40nL^2}\\
&= \frac{1}{8nL} - \frac{1}{40nL}\\
& \geqslant 0 \; .
\end{align*}
Hence,
\[
\E [ Q(\theta^k) | \mathcal{F}_{k-1} ]  - (1 - \delta)Q(\theta^{k-1}) \leqslant 0 \; .
\]

We can then take a full expectation on both sides to obtain:
\[
\E  Q(\theta^k)   - (1 - \delta) \E Q(\theta^{k-1}) \leqslant 0 \; .
\]

Since $Q$ is a non-negative function (we show below that it dominates a non-negative function), this results proves the linear convergence of the sequence $\E Q(\theta^k)$ with rate $1 - \delta$. We have
\begin{align*}
\E Q(\theta^k) 	& \leqslant \left(1 - \frac{\mu}{8nL}\right)^k Q(\theta^0) \; .
\end{align*}

\subsubsection*{Step 2 - Domination of $\|x^k - x^\ast\|^2$ by $Q(\theta^k)$}

We now need to prove that $Q(\theta^k)$ dominates $\|x^k - x^\ast\|^2$. If $\displaystyle P - \left( \begin{array}{cc} 0 & 0 \\ 0 & \frac{1}{3} \idm \end{array} \right)$ is positive definite, then $Q(\theta^k) \geqslant \frac{1}{3} \|x^k - x^\ast\|^2$.

We shall use the Schur complement condition for positive definiteness. Since $A$ is positive definite, the other condition to verify is $\frac{2}{3}\idm - b^\top A^{-1} b \succ 0$.
\begin{align*}
\frac{2}{3}\idm - \alpha^2\left(1 - \frac{1}{n}\right)^2 e^\top \left[\left(3 n\alpha^2 + \frac{\alpha^2}{n} - 2\alpha^2\right)\frac{ee^\top}{n}\right]^{-1}e &= \frac{2}{3}\idm -\frac{n\left(1 - \frac{1}{n}\right)^2}{3 n + \frac{1}{n} - 2}\frac{ee^\top}{n}\\
&\succ \frac{2}{3}\idm -\frac{n}{3 n - 2}\frac{ee^\top}{n}\\
&\succ 0 \textrm{ for } n \geqslant 2 \; ,
\end{align*}
and so $P$ dominates $\left( \begin{array}{cc} 0 & 0 \\ 0 & \frac{1}{3} \idm \end{array} \right)$.

This yields
\begin{align*}
\E \| x^k - x^\ast \|^2 	& \leqslant 3\E Q(\theta^k) \\
& \leqslant 3 \left(1 - \frac{\mu}{8nL}\right)^k Q(\theta^0) \; .
\end{align*}

We have
\begin{align*}
Q(\theta^0) &= 3n\alpha^2 \sum_i \|y_i^0 - f_i'(x^\ast)\|^2 + \frac{(1 - 2n)\alpha}{n^2} \left\|\sum_i y_i^0\right\|^2 - 2\alpha\left(1 - \frac{1}{n}\right)(x^0 - x^\ast)^\top\left(\sum_i y_i^0\right) + \|x^0 - x^\ast\|^2\\
&= \frac{3}{4nL^2}\sum_i \|y_i^0 - f_i'(x^\ast)\|^2 + \frac{(1 - 2n)}{2n^3L} \left\|\sum_i y_i^0\right\|^2 - \frac{n-1}{n^2L}(x^0 - x^\ast)^\top\left(\sum_i y_i^0\right) + \|x^0 - x^\ast\|^2 \; .
\end{align*}
Initializing all the $y_i^0$ to 0, we get
\begin{align*}
Q(\theta^0) &= \frac{3\sigma^2}{4L^2} + \|x^0 - x^\ast\|^2 \; ,
\end{align*}
and
\begin{align*}
\E \| x^k - x^\ast \|^2 & \leqslant \left(1 - \frac{\mu}{8nL}\right)^k \left(\frac{9\sigma^2}{4L^2} + 3\|x^0 - x^\ast\|^2\right) \; .
\end{align*}
\end{proof}

\subsection{Analysis for $\alpha = \frac{1}{2 n \mu}$}
\subsubsection*{Step 1 - Linear convergence of the Lyapunov function}
We now prove Proposition~\ref{prop:large_alpha}, providing a bound for the convergence rate of the SAG algorithm in the case of a small step size, $\alpha = \frac{1}{2n\mu}$.

We shall use the following Lyapunov function:
\[
Q(\theta^k) = 2 g \left(x^k + \frac{\alpha}{n} e^\top y^k\right) - 2 g(x^\ast) +
(\theta^k - \theta^\ast)^\top
\left(
\begin{array}{cc} A & b \\ b^\top & c
\end{array}
\right)
(\theta^k - \theta^\ast) \; ,
\]
with
\begin{align*}
A & = \frac{\eta \alpha}{n} \idm + \frac{\alpha}{n} (1 - 2 \nu )ee^\top \\
b & = - \nu e \\
c & = 0 \; .
\end{align*}

This yields
\begin{align*}
S & = \frac{ \eta \alpha}{n} \idm + \frac{\alpha}{n} ee^\top \\
\Diag(\diag(S)) & = \frac{ (1+\eta)\alpha}{n} \idm  \\
S - \Diag(\diag(S)) & = \frac{\alpha}{n} ( ee^\top - \idm )  \\
\left(1-\frac{2}{n}\right) S + \frac{1}{n} \Diag(\diag(S))
& = \left(1-\frac{2}{n}\right) \left[ \frac{ \eta \alpha}{n} \idm + \frac{\alpha}{n} ee^\top \right]
+ \frac{1}{n} \frac{ (1+\eta) \alpha}{n} \idm =
\left (1-\frac{2}{n}\right) \frac{\alpha}{n} ee^\top + \left(\eta - \frac{\eta - 1}{n}\right) \frac{\alpha}{n} \idm \; .
\end{align*}

We have
\begin{align*}
& \E [ Q(\theta^k) | \mathcal{F}_{k-1} ] - ( 1 - \delta) Q(\theta^{k-1}) \\
&\hspace*{.5cm} = 2 g (x^{k-1}) - 2 g(x^\ast) - 2 ( 1 - \delta) g \left(x^{k-1} + \frac{\alpha}{n} e^\top y^{k-1} \right)
+ 2 ( 1 - \delta) g( x^\ast) \\
&\hspace*{1cm} + (y^{k-1} - f'(x^\ast))^\top
\left[
\left(1-\frac{2}{n}\right) \frac{\alpha}{n} ee^\top + \left(\eta - \frac{\eta-1}{n}\right) \frac{\alpha}{n} \idm
- (1-\delta) \frac{\eta \alpha}{n} \idm \right. \\
& \hspace*{10cm} \left. - (1-\delta) \frac{\alpha}{n} (1 - 2 \nu )ee^\top
\right] (y^{k-1} - f'(x^\ast)) \\
&\hspace*{1cm} -\frac{2 \nu}{n} (x^{k-1} -x^\ast)^\top e ^\top ( f'(x^{k-1}) - f'(x^\ast)) \\
& \hspace*{1cm} + \frac{(1+\eta) \alpha}{n^2} ( f'(x^{k-1}) - f'(x^\ast))^\top  ( f'(x^{k-1}) - f'(x^\ast)) \\
&\hspace*{1cm} + \frac{2\alpha}{n^2} (y^{k-1} - f'(x^\ast))^\top \left[ ee^\top - \idm
\right] ( f'(x^{k-1}) - f'(x^\ast)) \\
&\hspace*{1cm} + 2 \left(\frac{1}{n}-\delta\right) \nu (y^{k-1} - f'(x^\ast))^\top e (x^{k-1} -x^\ast).
\end{align*}
Our goal will now be to express all the quantities in terms of $(x^{k-1} -x^\ast)^\top g'(x^{k-1})$ whose positivity is guaranteed by the convexity of $g$.

Using the convexity of $g$, we have
\begin{align*}
- 2 ( 1 - \delta) g \left(x^{k-1} + \frac{\alpha}{n} e^\top y^{k-1} \right) &\leqslant - 2 ( 1 - \delta) \left[ g(x^{k-1}) + \frac{\alpha}{n} g'(x^{k-1}) e^\top y^{k-1}
\right] \; .
\end{align*}
Using the Lipschitz property of the gradients of $f_i$, we have
\begin{align*}
( f'(x^{k-1}) - f'(x^\ast))^\top ( f'(x^{k-1}) - f'(x^\ast))
& = \sum_{i=1}^n \| f_i'(x^{k-1}) - f_i'(x^{\ast})\|^2 \\
&
\leqslant  \sum_{i=1}^n L ( f_i'(x^{k-1}) - f_i'(x^{\ast}))^\top( x^{k-1} - x^{\ast})
\\
& = nL ( g'(x^{k-1}) - g'(x^{\ast}))^\top( x^{k-1} - x^{\ast}) \; .
\end{align*}
Using $e^\top [ f'(x^{k-1}) - f'(x^\ast)] = ng'(x^{k-1})$, we have
\begin{align*}
-\frac{2 \nu}{n} (x^{k-1} -x^\ast)^\top e ^\top ( f'(x^{k-1}) - f'(x^\ast)) & = -2 \nu (x^{k-1} -x^\ast)^\top g'(x^{k-1})\\
\frac{2\alpha}{n^2} (y^{k-1} - f'(x^\ast))^\top ee^\top ( f'(x^{k-1}) - f'(x^\ast))&= \frac{2\alpha}{n} (y^{k-1} - f'(x^\ast))^\top e g'(x^{k-1}) \; .
\end{align*}
Reassembling all the terms together, we get
\begin{align*}
& \E [ Q(\theta^k) | \mathcal{F}_{k-1} ] - ( 1 - \delta) Q(\theta^{k-1}) \\
& \hspace*{.5cm}\leqslant
2 \delta[g (x^{k-1}) - g(x^\ast) ] + \frac{2\delta\alpha}{n} g'(x^{k-1}) e^\top y^{k-1}\\
&\hspace*{1cm} + (y^{k-1} - f'(x^\ast))^\top
\left[
\left(1-\frac{2}{n}\right) \frac{\alpha}{n} ee^\top + \left(\eta - \frac{\eta-1}{n}\right) \frac{\alpha}{n} \idm
- (1-\delta) \frac{\eta \alpha}{n} \idm - \right. \\
& \hspace*{10cm} \left. (1-\delta) \frac{\alpha}{n} (1 - 2 \nu )ee^\top
\right] (y^{k-1} - f'(x^\ast)) \\
&\hspace*{1cm} - \left( 2 \nu - \frac{(1+\eta) \alpha L}{n}\right)(x^{k-1} -x^\ast)^\top g'(x^{k-1}) \\
&\hspace*{1cm} - \frac{2\alpha}{n^2} (y^{k-1} - f'(x^\ast))^\top \big  ( f'(x^{k-1}) - f'(x^\ast)) \\
&\hspace*{1cm} + 2 \left(\frac{1}{n}-\delta\right) \nu (y^{k-1} - f'(x^\ast))^\top e (x^{k-1} -x^\ast).
\end{align*}
Using the convexity of $g$ gives
\begin{align*}
2 \delta[g (x^{k-1}) - g(x^\ast) ] &\leqslant 2 \delta[x^{k-1} -x^\ast ]^\top g'(x^{k-1}) \; ,
\end{align*}
and, consequently,
\begin{align*}
& \E [ Q(\theta^k) | \mathcal{F}_{k-1} ] - ( 1 - \delta) Q(\theta^{k-1}) \\
& \hspace*{.5cm}\leqslant
2 \delta[(x^{k-1}) -(x^\ast) ]^\top g'(x^{k-1}) + \frac{2\delta\alpha}{n} g'(x^{k-1}) e^\top y^{k-1}\\
&\hspace*{1cm} + (y^{k-1} - f'(x^\ast))^\top
\left[
\left(1-\frac{2}{n}\right) \frac{\alpha}{n} ee^\top + \left(\eta - \frac{\eta-1}{n}\right) \frac{\alpha}{n} \idm
\right. \\
& \hspace*{6cm} \left.
- (1-\delta) \frac{\eta \alpha}{n} \idm - (1-\delta) \frac{\alpha}{n} (1 - 2 \nu )ee^\top
\right] (y^{k-1} - f'(x^\ast)) \\
&\hspace*{1cm} - \left( 2 \nu - \frac{(1+\eta) \alpha L}{n}\right)(x^{k-1} -x^\ast)^\top g'(x^{k-1}) \\
&\hspace*{1cm} - \frac{2\alpha}{n^2} (y^{k-1} - f'(x^\ast))^\top \big  ( f'(x^{k-1}) - f'(x^\ast)) \\
&\hspace*{1cm} + 2 \left(\frac{1}{n}-\delta\right) \nu (y^{k-1} - f'(x^\ast))^\top e (x^{k-1} -x^\ast) \; .
\end{align*}
If we regroup all the terms in $[(x^{k-1}) -(x^\ast) ]^\top g'(x^{k-1})$ together, and all the terms in $(y^{k-1} - f'(x^\ast))^\top$ together, we get
\begin{align*}
& \E [ Q(\theta^k) | \mathcal{F}_{k-1} ] - ( 1 - \delta) Q(\theta^{k-1}) \\
& \hspace*{.5cm} \leqslant \frac{\alpha}{n} (y^{k-1} - f'(x^\ast))^\top
\left[\left( \delta \eta - \frac{\eta-1}{n} \right) \idm +\left( \delta -\frac{2}{n}  + 2 \nu (1-\delta) \right) ee^\top \right] (y^{k-1} - f'(x^\ast)) \\
&\hspace*{1cm} - \left( 2 \nu - 2 \delta - \frac{(1+\eta) \alpha L}{n} \right)(x^{k-1} -x^\ast)^\top g'(x^{k-1}) \\
&\hspace*{1cm} + 2 (y^{k-1} - f'(x^\ast))^\top\left[- \frac{\alpha}{n^2}  ( f'(x^{k-1}) - f'(x^\ast)) + (\frac{1}{n}-\delta) \nu e (x^{k-1} -x^\ast)
+\frac{\delta \alpha}{n} e g'(x^{k-1})  \right] \; .
\end{align*}

Let us rewrite this as
\begin{align*}
& \E [ Q(\theta^k) | \mathcal{F}_{k-1} ] - ( 1 - \delta) Q(\theta^{k-1}) \\
& \hspace*{1cm} \leqslant (y^{k-1} - f'(x^\ast))^\top \left(\tau_{y, I}I + \tau_{y, e}\frac{ee^\top}{n}\right)(y^{k-1} - f'(x^\ast))\\
& \hspace*{1cm}+ \tau_{x, g} (x^{k-1} -x^\ast)^\top g'(x^{k-1})\\
& \hspace*{1cm}+ (y^{k-1} - f'(x^\ast))^\top \left[\tau_{y, f}( f'(x^{k-1}) - f'(x^\ast)) + \tau_{y, x} e (x^{k-1} -x^\ast)+ \tau_{y, g} e g'(x^{k-1})\right]
\end{align*}
with
\begin{align*}
\tau_{y, I} &= \frac{\alpha}{n} \left(\delta\eta - \frac{\eta - 1}{n}\right)\\
\tau_{y, e} &= \alpha \left(\delta - \frac{2}{n} + 2 \nu(1 - \delta)\right)\\
\tau_{x,g} & = - ( 2 \nu - 2 \delta - \frac{(1+\eta) \alpha L}{n} )\\
\tau_{y,f} &= - \frac{2\alpha}{n^2}\\
\tau_{y,x}&= 2\left(\frac{1}{n}-\delta\right) \nu\\
\tau_{y,g}&= \frac{2\delta\alpha}{n}\; .
\end{align*}

Assuming that $\tau_{y, I}$ and $\tau_{y, e}$ are negative, we have by completing the square that
\begin{align*}
&(y^{k-1} - f'(x^\ast))^\top \left(\tau_{y, I}I + \tau_{y, e}\frac{ee^\top}{n}\right)(y^{k-1} - f'(x^\ast))\\
&\hspace*{1cm} + (y^{k-1} - f'(x^\ast))^\top\left(\tau_{y, f}( f'(x^{k-1}) - f'(x^\ast)) + \tau_{y, x}e (x^{k-1} -x^\ast) + \tau_{y, g}e g'(x^{k-1})\right)\\
&\hspace*{.5cm}\leqslant -\frac{1}{4} \left(\tau_{y, f}( f'(x^{k-1}) - f'(x^\ast)) + \tau_{y, x}e (x^{k-1} -x^\ast) + \tau_{y, g}e g'(x^{k-1})\right)^\top\left(\frac{1}{\tau_{y,I}}\left(I - \frac{ee^\top}{n}\right) + \frac{1}{\tau_{y,I} + \tau_{y,e}}\frac{ee^\top}{n}\right)\\
&\hspace*{1cm}\left(\tau_{y, f}( f'(x^{k-1}) - f'(x^\ast)) + \tau_{y, x}e (x^{k-1} -x^\ast) + \tau_{y, g}e g'(x^{k-1})\right)\\
&\hspace*{.5cm}= -\frac{1}{4}\frac{\tau_{y, f}^2}{\tau_{y,I}}\|f'(x^{k-1}) - f'(x^\ast)\|^2 -\frac{1}{4}\tau_{y, f}^2n \|g'(x^{k-1})\|^2\left(\frac{1}{\tau_{y,I} + \tau_{y,e}} - \frac{1}{\tau_{y,I}}\right)\\
&\hspace*{1cm}-\frac{1}{4}\frac{\tau_{y,x}^2n}{\tau_{y,I} + \tau_{y,e}}\|x^{k-1} -x^\ast\|^2 -\frac{1}{4}\frac{\tau_{y,g}^2n}{\tau_{y,I} + \tau_{y,e}}\|g'(x^{k-1})\|^2\\
&\hspace*{1cm}-\frac{1}{2}\frac{\tau_{y,f}\tau_{y,x}n}{\tau_{y,I} + \tau_{y,e}}(x^{k-1} -x^\ast)^\top g'(x^{k-1}) -\frac{1}{2}\frac{\tau_{y,f}\tau_{y,g}n}{\tau_{y,I} + \tau_{y,e}}\|g'(x^{k-1}) \|^2 -\frac{1}{2}\frac{\tau_{y,g}\tau_{y,x}n}{\tau_{y,I} + \tau_{y,e}}(x^{k-1} -x^\ast)^\top g'(x^{k-1}) \; ,
\end{align*}
where we used the fact that $( f'(x^{k-1}) - f'(x^\ast))^\top e = g'(x^{k-1})$. After reorganization of the terms, we obtain
\begin{align*}
\E [ Q(\theta^k) | \mathcal{F}_{k-1} ] - ( 1 - \delta) Q(\theta^{k-1}) &\leqslant \left[\tau_{x,g} - \frac{n\tau_{y,x}}{2(\tau_{y,I} + \tau_{y,e})}(\tau_{y,f} + \tau_{y,g}) \right](x^{k-1} -x^\ast)^\top g'(x^{k-1})\\
& \hspace*{.5cm} -\left[\frac{1}{4}\tau_{y, f}^2n\left(\frac{1}{\tau_{y,I} + \tau_{y,e}} - \frac{1}{\tau_{y,I}}\right) + \frac{1}{4}\frac{\tau_{y,g}^2n}{\tau_{y,I} + \tau_{y,e}} + \frac{1}{2}\frac{\tau_{y,f}\tau_{y,g}n}{\tau_{y,I} + \tau_{y,e}}\right]\|g'(x^{k-1}) \|^2\\
& \hspace*{.5cm} -\frac{1}{4}\frac{\tau_{y, f}^2}{\tau_{y,I}}\|f'(x^{k-1}) - f'(x^\ast)\|^2 -\frac{1}{4}\frac{\tau_{y,x}^2n}{\tau_{y,I} + \tau_{y,e}}\|x^{k-1} -x^\ast\|^2 \; .
\end{align*}

We now use the strong convexity of the function to get the following inequalities:
\begin{align*}
\|f'(x^{k-1}) - f'(x^\ast)\|^2 & \leqslant Ln (x^{k-1} -x^\ast)^\top g'(x^{k-1})\\
\|x^{k-1} -x^\ast\|^2& \leqslant\frac{1}{\mu}(x^{k-1} -x^\ast)^\top g'(x^{k-1}) \; .
\end{align*}

Finally, we have
\begin{align*}
& \E [ Q(\theta^k) | \mathcal{F}_{k-1} ] - ( 1 - \delta) Q(\theta^{k-1}) \\
& \hspace*{.5cm} \leqslant \left[\tau_{x,g} - \frac{n\tau_{y,x}}{2(\tau_{y,I} + \tau_{y,e})}(\tau_{y,f} + \tau_{y,g}) -\frac{Ln}{4}\frac{\tau_{y, f}^2}{\tau_{y,I}}-\frac{1}{4\mu}\frac{\tau_{y,x}^2n}{\tau_{y,I} + \tau_{y,e}}\right](x^{k-1} -x^\ast)^\top g'(x^{k-1})\\
& \hspace*{1cm} -\left[\frac{1}{4}\tau_{y, f}^2n\left(\frac{1}{\tau_{y,I} + \tau_{y,e}} - \frac{1}{\tau_{y,I}}\right) + \frac{1}{4}\frac{\tau_{y,g}^2n}{\tau_{y,I} + \tau_{y,e}} + \frac{1}{2}\frac{\tau_{y,f}\tau_{y,g}n}{\tau_{y,I} + \tau_{y,e}}\right]\|g'(x^{k-1}) \|^2 \; .
\end{align*}

If we choose $\delta = \frac{\widetilde{\delta}}{n}$ with $\widetilde{\delta} \leqslant \frac{1}{2}$, $\nu = \frac{1}{2n}$, $\eta = 2$ and $\alpha = \frac{1}{2n\mu}$, we get
\begin{align*}
\tau_{y, I} &= \frac{1}{2n^2\mu} \left(\frac{2\widetilde{\delta}}{n} - \frac{1}{n}\right) = -\frac{1 - 2\widetilde{\delta}}{2n^3\mu} \leqslant 0\\
\tau_{y, e} &=\frac{1}{2n\mu} \left( \frac{\widetilde{\delta}}{n} - \frac{2}{n} + \frac{1}{n}\left(1 - \frac{\widetilde{\delta}}{n}\right)\right) = -\frac{1}{2n^2\mu}\left(1 - \widetilde{\delta} + \frac{\widetilde{\delta}}{n}\right) \leqslant 0\\
\tau_{x,g} & = -\left ( \frac{1}{n} - \frac{2\widetilde{\delta}}{n} - \frac{3L}{2n^2\mu} \right) = \frac{3 L}{2n^2\mu} - \frac{1 - 2\widetilde{\delta}}{n}\\
\tau_{y,f} &= - \frac{1}{n^3\mu}\\
\tau_{y,x}&= \frac{1 - \widetilde{\delta}}{n^2}\\
\tau_{y,g}&= \frac{\widetilde{\delta}}{n^3\mu}\; .
\end{align*}

Thus,
\begin{align*}
\tau_{x,g}&- \frac{n\tau_{y,x}}{2(\tau_{y,I} + \tau_{y,e})}(\tau_{y,f} + \tau_{y,g}) -\frac{Ln}{4}\frac{\tau_{y, f}^2}{\tau_{y,I}}-\frac{1}{4\mu}\frac{\tau_{y,x}^2n}{\tau_{y,I} + \tau_{y,e}}\\
&\leqslant \frac{3 L}{2n^2\mu} - \frac{1 - 2\widetilde{\delta}}{n} - \frac{\frac{1 - \widetilde{\delta}}{2n}\frac{2\widetilde{\delta} - 1}{n^3\mu}}{\tau_{y,I} + \tau_{y,e}} + \frac{Ln}{4}\frac{\frac{1}{n^6\mu^2}}{\frac{1 - 2\widetilde{\delta}}{2n^3\mu}} - \frac{1}{4\mu}\frac{\frac{(1 - \widetilde{\delta})^2}{n^3}}{\tau_{y,I} + \tau_{y,e}}\\
& = \frac{L}{n^2\mu}\left[\frac{3}{2} + \frac{1}{2(1 - 2\widetilde{\delta})}\right] - \frac{1 - 2\widetilde{\delta}}{n}- \frac{1}{\mu n^3 (\tau_{y,I} + \tau_{y,e})} \left[\frac{(1-\widetilde{\delta})^2}{4} + \frac{(1 - \widetilde{\delta})(2\widetilde{\delta} - 1)}{2n}\right]\\
& \leqslant \frac{L}{n^2\mu} \frac{2 - 3\widetilde{\delta}}{1 - 2\widetilde{\delta}} - \frac{1 - 2\widetilde{\delta}}{n}+ \frac{1}{\mu n^3 \left(\frac{1 - 2\widetilde{\delta}}{2n^3\mu} + \frac{1}{2n^2\mu}\left(1 - \widetilde{\delta} + \frac{\widetilde{\delta}}{n}\right)\right)}\frac{(1-\widetilde{\delta})^2}{4}\\
&= \frac{L}{n^2\mu} \frac{2 - 3\widetilde{\delta}}{1 - 2\widetilde{\delta}} - \frac{1 - 2\widetilde{\delta}}{n} + \frac{(1 - \widetilde{\delta})^2}{2 - 4\widetilde{\delta} + 2n - 2n\widetilde{\delta} + 2\widetilde{\delta}}\\
&= \frac{L}{n^2\mu} \frac{2 - 3\widetilde{\delta}}{1 - 2\widetilde{\delta}} - \frac{1 - 2\widetilde{\delta}}{n} + \frac{1 - \widetilde{\delta}}{2(1 + n)}\\
&\leqslant \frac{L}{n^2\mu} \frac{1 - 3\widetilde{\delta}}{1 - 2\widetilde{\delta}} - \frac{1 - 2\widetilde{\delta}}{n} + \frac{1 - \widetilde{\delta}}{2n}\\
& = \frac{L}{n^2\mu} \frac{2 - 3\widetilde{\delta}}{1 - 2\widetilde{\delta}} - \frac{1- 3\widetilde{\delta}}{2n} \; .
\end{align*}
This quantity is negative for $\widetilde{\delta} \leqslant \frac{1}{3}$ and $\frac{\mu}{L} \geqslant \frac{4 - 6\widetilde{\delta}}{n(1 - 2\widetilde{\delta})(1 - 3\widetilde{\delta})}$.
If we choose $\widetilde{\delta} = \frac{1}{8}$, then it is sufficient to have $\frac{n\mu}{L} \geqslant 8$.

To finish the proof, we need to prove the positivity of the factor of $\|g'(x^{k-1}) \|^2$.
\begin{align*}
\frac{1}{4}\tau_{y, f}^2n\left(\frac{1}{\tau_{y,I} + \tau_{y,e}} - \frac{1}{\tau_{y,I}}\right) + \frac{1}{4}\frac{\tau_{y,g}^2n}{\tau_{y,I} + \tau_{y,e}} + \frac{1}{2}\frac{\tau_{y,f}\tau_{y,g}n}{\tau_{y,I} + \tau_{y,e}}& = \frac{n}{4}\frac{1}{\tau_{y,I} + \tau_{y,e}}(\tau_{y, f} + \tau_{y,g})^2 - \frac{n}{4}\frac{\tau_{y, f}^2}{\tau_{y,I}}\\
&\geqslant \frac{n}{4}\frac{(\tau_{y, f} + \tau_{y,g})^2}{\tau_{y,I}} - \frac{n}{4}\frac{\tau_{y, f}^2}{\tau_{y,I}}\\
&= \frac{n}{4\tau_{y,I}} \tau_{y,g}(2\tau_{y, f} + \tau_{y,g})\\
&\geqslant 0 \; .
\end{align*}

Then, following the same argument as in the previous section, we have
\begin{align*}
\E Q(\theta^k) 	&\leqslant \left(1 - \frac{1}{8n}\right)^k Q(\theta^0)\\
& = \left(1 - \frac{1}{8n}\right)^k \left[2(g(x^0) - g(x^\ast)) + \frac{\sigma^2}{n\mu}\right] \; ,
\end{align*}
with $\sigma^2 = \frac{1}{n} \sum_i \|f_i'(x^\ast)\|^2$ the variance of the gradients at the optimum.

\subsubsection*{Step 2 - Domination of $g(x^k) - g(x^\ast)$ by $Q(\theta^k)$}

We now need to prove that $Q(\theta^k)$ dominates $g(x^k) - g(x^\ast)$.

\begin{align*}
Q(\theta^k) &= 2 g \left(x^k + \frac{\alpha}{n} e^\top y^k\right) - 2 g(x^\ast) +
(\theta^k - \theta^\ast)^\top
\left(
\begin{array}{cc} A & b \\ b^\top & c
\end{array}
\right)
(\theta^k - \theta^\ast)\\
&= 2 g \left(x^k + \frac{\alpha}{n} e^\top y^k\right) - 2 g(x^\ast) + \frac{1}{n^2\mu}\sum_i\left\|y_i^k - f_i'(x^\ast)\right\|^2 + \frac{n-1}{2n^3\mu}\|e^\top y\|^2 - \frac{1}{n} (x^k - x^\ast)^\top (e^\top y^k)\\
& \geqslant 2 g(x^k) + \frac{2\alpha}{n}g'(x^k)^\top (e^\top y^k) - 2 g(x^\ast)\\
&\hspace*{.5cm} + \frac{1}{n^2\mu}\sum_i\left\|\frac{1}{n}e^\top y^k + y_i^k -\frac{1}{n}e^\top y^k - f_i'(x^\ast)\right\|^2 + \frac{n-1}{2n^3\mu}\|e^\top y\|^2 - \frac{1}{n} (x^k - x^\ast)^\top (e^\top y^k)\\
& \mbox{ using the convexity of } g \mbox{ and the fact that } \sum_{i} f_i'(x^\ast)=0\\
&= 2 g(x^k) - 2 g(x^\ast) + \left(\frac{2\alpha}{n}g'(x^k) - \frac{1}{n} (x^k - x^\ast)\right)^\top (e^\top y^k)\\
&\hspace*{.5cm} + \frac{1}{n^3\mu} \|e^\top y^k\|^2 + \frac{1}{n^2\mu}\sum_i\left\|y_i^k -\frac{1}{n}e^\top y^k - f_i'(x^\ast)\right\|^2 + \frac{n-1}{2n^3\mu}\|e^\top y\|^2\\
&\geqslant 2 g(x^k) - 2 g(x^\ast) + \left(\frac{2\alpha}{n}g'(x^k) - \frac{1}{n} (x^k - x^\ast)\right)^\top (e^\top y^k) + \frac{n+1}{2n^3\mu}\|e^\top y\|^2 \; \\
& \mbox{ by dropping some terms.}\\
\end{align*}

The quantity on the right-hand side is minimized for $e^\top y = \frac{n^3\mu}{n+1}\left(\frac{1}{n} (x^k - x^\ast) - \frac{2\alpha}{n}g'(x^k)\right)$. Hence, we have

\begin{align*}
Q(\theta^k) &\geqslant 2 g(x^k) - 2 g(x^\ast) - \frac{n^3\mu}{2(n+1)} \left\|\frac{1}{n} (x^k - x^\ast) - \frac{2\alpha}{n}g'(x^k)\right\|^2\\
&= 2 g(x^k) - 2 g(x^\ast) - \frac{n^3\mu}{2(n+1)} \left( \frac{1}{n^2} \|x^k - x^\ast\|^2 + \frac{4\alpha^2}{n^2}\|g'(x^k)\|^2 - \frac{4\alpha}{n^2} (x^k - x^\ast)^\top g'(x^k)\right)\\
& \geqslant 2 g(x^k) - 2 g(x^\ast) - \frac{n^3\mu}{2(n+1)} \left( \frac{1}{n^2} \|x^k - x^\ast\|^2 + \frac{4\alpha^2}{n^2}\|g'(x^k)\|^2\right)\\
& \mbox{ using the convexity of } g \\
&\geqslant 2 g(x^k) - 2 g(x^\ast) - \frac{n\mu}{2(n+1)} \left(1 + \frac{L^2}{\mu^2n^2}\right)\|x^k - x^\ast\|^2\\
& \mbox{ using the Lipschitz continuity of } g' \\
&\geqslant 2 g(x^k) - 2 g(x^\ast) - \frac{n\mu}{2(n+1)} \frac{65}{64}\|x^k - x^\ast\|^2 \textrm{ since } \frac{\mu}{L} \geqslant \frac{8}{n}\\
&\geqslant 2 g(x^k) - 2 g(x^\ast) - \frac{n}{(n+1)} \frac{65}{64}(g(x^k) - g(x^\ast))\\
&\geqslant \frac{63}{64}(g(x^k) - g(x^\ast))\\
&\geqslant \frac{6}{7}(g(x^k) - g(x^\ast)) \; .
\end{align*}

We thus get
\begin{align*}
\E\left[g(x^k) - g(x^\ast)\right] 	& \leqslant 2 \E Q(\theta^k)\\
& = \left(1 - \frac{1}{8n}\right)^k \left[\frac{7}{3}(g(x^0) - g(x^\ast)) + \frac{7\sigma^2}{6n\mu}\right] \; .
\end{align*}

\subsubsection*{Step 3 - Initialization of $x^0$ using stochastic gradient descent}
During the first few iterations, we obtain the $O(1/k)$ rate obtained using stochastic gradient descent, but with a constant which is proportional to $n$. To circumvent this problem, we will first do $n$ iterations of stochastic gradient descent to initialize $x^0$, which will be renamed $x^n$ to truly reflect the number of iterations done.

Using the bound from section~\ref{sec:bound_values_sg}, we have
\begin{align*}
\E g\left(\frac{1}{n}\sum_{i=0}^{n-1}\tilde{x}^{i}\right) - g(x^\ast) & \leqslant \frac{2L}{n}\| x^0 - x^\ast\|^2 + \frac{4 \sigma^2 }{n \mu} \log\left( 1 + \frac{\mu n }{4L } \right) \; .
\end{align*}
And so, using $x^n = \frac{1}{n}\sum_{i=0}^{n-1}\tilde{x}^{i}$, we have for $k \geqslant n$
\begin{align*}
\E\left[g(x^k) - g(x^\ast)\right] & \leqslant \left(1 - \frac{1}{8n}\right)^{k-n} \left[\frac{14L}{3n}\| x^0 - x^\ast\|^2 + \frac{28 \sigma^2 }{3n \mu} \log\left( 1 + \frac{\mu n }{4L } \right) + \frac{7\sigma^2}{6n\mu}\right] \; .
\end{align*}
Since
\begin{align*}
\left(1 - \frac{1}{8n}\right)^{-n} &\leqslant \frac{8}{7} \; ,
\end{align*}
we get
\begin{align*}
\E\left[g(x^k) - g(x^\ast)\right] & \leqslant \left(1 - \frac{1}{8n}\right)^k \left[\frac{16L}{3n}\| x^0 - x^\ast\|^2 + \frac{32 \sigma^2 }{3n \mu} \log\left( 1 + \frac{\mu n }{4L } \right) + \frac{4\sigma^2}{3n\mu}\right] \; .
\end{align*}

\section{Comparison of convergence rates}

We consider the $\ell_2$-regularized least squares problem
\[
\minimize{x\in\Real^p} \quad g(x) \defd \frac{\lambda}{2}\norm{x}^2 + \frac{1}{2n}\sum_{i=1}^n(a_i^Tx - b_i)^2,
\]
where to apply SG methods and SAG we can use
\[
f_i(x) := \frac{\lambda}{2}\norm{x}^2 + \frac{1}{2}(a_i^Tx - b_i)^2.
\]
If we use $b$ to denote a vector containing the values $b_i$ and $A$ to denote a matrix withs rows $a_i$, we can re-write this problem as
\[
\minimize{x\in\Real^p} \quad \frac{\lambda}{2}\norm{x}^2 + \frac{1}{2n}\norm{Ax - b}^2.
\]
The Fenchel dual of this problem is
\[
\minimize{y\in\Real^n} \quad d(y) \defd \frac{n}{2}\norm{y}^2 + \frac{1}{2\lambda}y^\top AA^\top y + y^\top b.
\]
We can obtain the primal variables from the dual variables by the formula $x = (-1/\lambda)A^\top y$.
Convergence rates of different primal and dual algorithms are often expressed in terms of the following Lipschitz constants:
\begin{align*}
&L_g = \lambda + M_\sigma/n & & \textrm{(Lipschitz constant of $g'$)}\\
&L_g^i = \lambda + M_i \quad & & \textrm{(Lipschitz constant for all $f_i'$)}\\
&L_g^j = \lambda + M_j/n \quad & & \textrm{(Lipschitz constant of all $g_j'$)}\\
&L_d = n + M_\sigma/\lambda & & \textrm{(Lipschitz constant of $d'$)}\\
&L_d^i = n + M_i/\lambda \quad & & \textrm{(Lipschitz constant of all $d_i'$)}\\
\end{align*}
Here, we use $M_\sigma$ to denote the maximum eigenvalue of $A^\top A$, $M_i$ to denote the maximum squared row-norm $\max_i\{\norm{a_i}^2\}$, and $M_j$ to denote the maximum squared column-norm $\max_j\{\sum_{i=1}^n(a_i)^2_j\}$. We use $g_j'$ to refer to element of $j$ of $g'$, and similarly for $d_i'$. The convergence rates will also depend on the primal and dual strong-convexity constants:
\begin{align*}
&\mu_g = \lambda + m_\sigma/n & & \textrm{(Strong-convexity constant of $g$)}\\
&\mu_d = n + m_\sigma'/\lambda & & \textrm{(Strong-convexity constant of $d$)}
\end{align*}
Here, $m_\sigma$ is the minimum eigenvalue of $A^\top A$, and $m_\sigma'$ is the minimum eigenvalue of $AA^\top$.

\subsection{Full Gradient Methods}

Using a similar argument to~\cite[Theorem 2.1.15]{nesterov2004introductory}, if we use the basic FG method with a step size of $1/L_g$, then $(f(x^k) - f(x^\ast))$ converges to zero with rate
\[
\left(1-\frac{\mu_g}{L_g}\right)^2 
= \left(1 - \frac{\lambda + m_\sigma/n}{\lambda + M_\sigma/n}\right)^2 
= \left(1 - \frac{n\lambda + m_\sigma}{n\lambda + M_\sigma}\right)^2 
\leq \exp\left(-2\frac{n\lambda + m_\sigma}{n\lambda + M_\sigma}\right),
\]
while a larger step-size of $2/(L_g+\mu_g)$ gives a faster rate of
\[
\left(1 - \frac{\mu_g + \mu_g}{L_g + \mu_g}\right)^2
= \left(1 - \frac{n\lambda + m_\sigma}{n\lambda + (M_\sigma + m_\sigma)/2}\right)^2
\leq \exp\left(-2\frac{n\lambda _+ m_\sigma}{n\lambda + (M_\sigma + m_\sigma)/2}\right),
\]
where the speed improvement is determined by the size of $m_\sigma$.

If we use the basic FG method on the dual problem with a step size of $1/L_d$, then $(d(x^k) - d(x^\ast))$ converges to zero with rate
\[
\left(1-\frac{\mu_d}{L_d}\right)^2
= \left(1-\frac{n + m_\sigma'/\lambda}{n + M_\sigma/\lambda}\right)^2 
= \left(1-\frac{n\lambda + m_\sigma'}{n\lambda + M_\sigma}\right)^2
\leq \exp\left(-2\frac{n\lambda + m'_\sigma}{n\lambda + M_\sigma}\right),
\]
and with a step-size of $2/(L_d+\mu_d)$ the rate is
\[
\left(1 - \frac{\mu_d + \mu_d}{L_d + \mu_d}\right)^2
= \left(1 - \frac{ n\lambda + m_\sigma'}{n\lambda + (M_\sigma + m_\sigma')/2}\right)^2
\leq \exp\left(-2\frac{n\lambda + m'_\sigma}{n\lambda + (M_\sigma + m'_\sigma)/2}\right).
\]
Thus, whether we can solve the primal or dual method faster depends on $m_\sigma$ and $m_\sigma'$. In the over-determined case where $A$ has independent columns, a primal method should be preferred. In the under-determined case where $A$ has independent rows, we can solve the dual more efficiently. However, we note that a convergence rate on the dual objective does not necessarily yield the same rate in the primal objective. If $A$ is invertible (so that $m_\sigma=m_\sigma'$) or it has neither independent columns nor independent rows (so that $m_\sigma = m_\sigma' = 0$), then there is no difference between the primal and dual rates.

The AFG method achieves a faster rate. Applied to the primal with a step-size of $1/L_g$ it has a rate of~\cite[Theorem 2.2.2]{nesterov2004introductory}
\[
\left(1-\sqrt{\frac{\mu_g}{L_g}}\right)
= \left(1-\sqrt{\frac{\lambda + m_\sigma/n}{\lambda + M_\sigma/n}}\right)
= \left(1-\sqrt{\frac{n\lambda + m_\sigma}{n\lambda + M_\sigma}}\right)
\leq \exp\left(-\sqrt{\frac{n\lambda + m_\sigma}{n\lambda + M_\sigma}}\right),
\]
and applied to the dual with a step-size of $1/L_d$ it has a rate of
\[
\left(1-\sqrt{\frac{\mu_d}{L_d}}\right)
= \left(1-\sqrt{\frac{n + m_\sigma'\lambda}{n + M_\sigma/\lambda}}\right)
= \left(1-\sqrt{\frac{n\lambda + m_\sigma'}{n\lambda + M_\sigma}}\right)
\leq \exp\left(-\sqrt{\frac{n\lambda + m'_\sigma}{n\lambda + M_\sigma}}\right).
\]
\subsection{Coordinate-Descent Methods}

The cost of applying one iteration of an FG method is $O(np)$. For this same cost we could apply $p$ iterations of a coordinate descent method to the primal, assuming that selecting the coordinate to update has a cost of $O(1)$. If we select coordinates uniformly at random, then the convergence rate for $p$ iterations of coordinate descent with a step-size of $1/L_g^j$ is~\cite[Theorem 2]{nesterov2010efficiency}
\[
\left(1-\frac{\mu_g}{pL_g^j}\right)^p
= \left(1-\frac{\lambda + m_\sigma/n}{p(\lambda + M_j/n)}\right)^p
= \left(1-\frac{n\lambda + m_\sigma}{p(n\lambda + M_j)}\right)^p
\leq \exp\left(-\frac{n\lambda + m_\sigma}{n\lambda + M_j}\right).
\]
Here, we see that applying a coordinate-descent method can be much more efficient than an FG method if $M_j << M_\sigma$. This can happen, for example, when the number of variables $p$ is much larger than the number of examples $n$. Further, it is possible for coordinate descent to be faster than the AFG method if the difference between $M_\sigma$ and $M_j$ is sufficiently large.

For the $O(np)$ cost of one iteration of the FG method, we could alternately perform $n$ iterations of coordinate descent on the dual problem. With a step size of $1/L_d^i$ this would obtain a rate on the dual objective of
\[
\left(1-\frac{\mu_d}{nL_d^i}\right)^n
= \left(1-\frac{n + m_\sigma'/\lambda}{n(n + M_i/\lambda)}\right)^n
= \left(1-\frac{n\lambda + m_\sigma'}{n(n\lambda + M_i)}\right)^n
\leq \exp\left(-\frac{n\lambda + m'_\sigma}{n\lambda + M_i}\right),
\]
which will be faster than the dual FG method if $M_i << M_\sigma$. This can happen, for example, when the number of examples $n$ is much larger than the number of variables $p$. The difference between the primal and dual coordinate methods depends on $M_i$ compared to $M_j$ and $m_\sigma$ compared to $m_\sigma'$.

\subsection{Stochastic Average Gradient}

For the $O(np)$ cost of one iteration of the FG method, we can perform $n$ iterations of SAG. With a step size of $1/2nL_g$, performing $n$ iterations of the SAG algorithm has a rate of
\[
\left(1 - \frac{\mu_g}{8nL_g^i}\right)^n
= \left(1 - \frac{\lambda + m_\sigma/n}{8n(\lambda + M_i)}\right)^n
= \left(1 - \frac{n\lambda + m_\sigma}{8n(n\lambda + nM_i)}\right)^n
\leq \exp\left(-\frac{1}{8}\frac{n\lambda + m_\sigma}{n\lambda + nM_i}\right),
\]
This is most similar to the rate obtained with the dual coordinate descent method, but is likely to be slower because of the $n$ term scaling $M_i$. However, the difference will be decreased for over-determined problems when $m_\sigma >> m_\sigma'$.

Under the condition $n \geqslant 8L_g^i/\mu_g = 8(\lambda + M_i)/(\lambda + m_\sigma/n)$, with a step size of $1/2n\mu_g$ performing $n$ iterations of the SAG algorithm has a rate of
\[
\left(1 - \frac{1}{8n}\right)^n
= \left(1 - \frac{n\lambda}{8n(n\lambda)}\right)^n
\leq \exp\left(-\frac{1}{8}\right).
\]
Note that depending on the constants this may or may not not be faster than coordinate descent methods.  However, if we consider the typical case where $m_\sigma = m_\sigma' = 0$ with $M_i = O(p)$ and $M_j = O(n)$, then if we have $n = 8(\lambda + M_i)/\lambda$ we obtain
\[
\left(1 - \frac{1}{8n}\right)^n
= \left(1 - \frac{\lambda}{64(\lambda + M_i)}\right)^n
= \left(1 - \frac{n\lambda}{64n(\lambda + M_i)}\right)^n
\leq \exp\left(-\frac{1}{64}\frac{n\lambda}{\lambda + M_i}\right).
\]
Despite the constant of $64$ (which is likely to be highly sub-optimal), from these rates we see that SAG outperforms coordinate descent methods when $n$ is sufficiently large.

\bibliographystyle{abbrvnat}

\setlength{\parsep}{0pt}
\setlength{\parskip}{0pt}
\bibliography{bibnips}

\end{document}